\RequirePackage{fix-cm}
\documentclass{svjour3}                     
\smartqed  

\usepackage{listings}
\usepackage{pdflscape}
\usepackage{epsfig}
\usepackage{lscape}
\usepackage{longtable}
\usepackage{amsmath,amssymb}
\usepackage{mathabx}
\usepackage{graphicx}
\usepackage{color}
\usepackage{placeins}
\usepackage{url}
\usepackage{enumitem,multirow}
\usepackage{framed}
\usepackage[figuresright]{rotating}
\usepackage{cite}
\usepackage{url}
\usepackage{hyperref}
\usepackage{algorithm}
\usepackage{algorithmic}
\usepackage{subfigure}
\usepackage[justification=centering]{caption}
\newtheorem{assumption}{Assumption}

\def\nn{\nonumber}

\DeclareMathOperator*{\argmin}{argmin}

\graphicspath{ {./figures/} }

\begin{document}

\title{A proximal DC approach for quadratic assignment problem}
\author{Zhuoxuan Jiang         \and
        Xinyuan Zhao \and
        Chao Ding
}
\institute{Z.X. Jiang \at
              College of Applied Sciences, Beijing University of Technology, Beijing, P.R. China. \\
              \email{zxjiang@emails.bjut.edu.cn}           
                         \and
           X.Y. Zhao \at
              College of Applied Sciences, Beijing University of Technology, Beijing, P.R. China. The research of this author was supported by the National Natural Science Foundation of China under projects No. 11871002 and the General Program of Science and Technology of Beijing Municipal Education Commission.\\ 
              \email{xyzhao@bjut.edu.cn}
              \and
              C. Ding \at
              Institute of Applied Mathematics, Academy of Mathematics and Systems Science, Chinese Academy of Sciences, Beijing,  P.R. China. The research of this author was supported by the National Natural Science Foundation of China under projects
		No. 11671387 and No. 11531014.\\ 
              \email{dingchao@amss.ac.cn}
}

\date{August 13, 2019}

\maketitle

\begin{abstract}
In this paper,  we show that the quadratic assignment problem (QAP) can be reformulated to an equivalent rank constrained doubly nonnegative (DNN) problem. Under the framework of the difference of convex functions (DC) approach, a semi-proximal DC algorithm (DCA) is proposed for solving the relaxation of the rank constrained DNN problem whose subproblems can be solved by the semi-proximal augmented Lagrangian method (sPALM). We show that the generated sequence converges to a stationary point of the corresponding DC problem, which is feasible to the rank constrained DNN problem. Moreover,  numerical experiments demonstrate that for most QAP instances, the proposed approach can find the global optimal solutions efficiently, and for others, the proposed algorithm is able to provide good feasible solutions in a reasonable time.  

\end{abstract}

\section{Introduction}

The quadratic assignment problem (QAP) is a classical mathematical model for location theory, which is used to model the location problem of allocating $n$ facilities to $n$ locations while minimizing the quadratic objective coming from the distance between the locations and the flow between the facilities. The standard form introduced by Koopmans and Beckmann \cite{KBeckmann57} is as following:
\begin{equation}\label{eq:QAP-original}
	\min \left\{ \sum_{1\le i,j\le n}A_{ij}B_{\pi(i),\pi(j)}+\sum_{i}C_{i\pi(i)}\mid \pi\in {\cal P}^n \right\},
\end{equation}
where $A$, $B$ and $C$ are given $n\times n$ real matrices and ${\cal P}^n$ is the the group of all permutations of $\{ 1, \ldots, n\}$. In this paper, we make the standard assumption that $A$ and $B$ are symmetric. 

Nowadays, QAP becomes one of the most important combinatorial optimization problems due to its widely applications in many different areas, such as chip design, manufacturing, computer graphics and vision, and so on (see \cite{Burkard13,Drezner15} for more details). However, it is well known that QAP is NP-hard \cite{SGonzalez76} and still quite difficult to compute the problems of dimension $n \ge 30$ in a reasonable computational time. Exact solution algorithms for QAP in practice are usually based on the branch and bound technique which is used to reduce the domain and to improve the bounds of relaxation problems \cite{Anstreicher03}.  Therefore, it is still an important research topic to improve the lower or upper bounds for QAP efficiently.

Meanwhile, semidefinite programming (SDP) \cite{Todd01} has proven to be very successful in this trend by providing tight relaxations for hard combinatorial problems\cite{VBoyd96}. To obtain lower bounds for QAP, various SDP relaxations are established \cite{LSaigal97,ZhaoKRendlW98}.  Although SDP relaxation is numerically successful,  it does not satisfy the Slater condition that may make the dual optimal solution unbounded \cite{RTWolkowicz97}. That is an important reason why some interior-point methods become inefficient for solving QAPs. To overcome this difficulty, by exploring the geometrical structure of SDP relaxations, Zhao et al. \cite{ZhaoKRendlW98} considered a reduced SDP problem by projecting the primal problem onto the minimal face of the semidefinite cone, and  constructed some Slater points for such SDP relaxations, which can be solved by the interior-point method and the bundle method \cite{RSot07} efficiently for $n \le 30$. 

In order to improve the quality of the SDP relaxation of QAP, Povh and Rendl \cite{PRendl09} showed that the optimal value of QAP was equal to the optimal value of the convex completely positive programming (CPP), i.e., a linear program over the cone of completely positive matrices. In fact, based on \cite{Bur09}, many important binary and nonconvex quadratic programs including QAP can be equivalent reformulated as the convex CPPs, under some mild conditions. However, these CPP reformulations are known to be numerically intractable \cite{MKab87}, and an efficient strategy is replacing the completely positive cone with doubly nonnegative (DNN) cone and solving the relaxation problems by SDP solvers \cite{FGYe18, YMat10,ZhaoSunToh10, WGYin10, KKojimaToh15, YangSunToh15}.
The QAP and the corresponding CPP relaxation proposed by Povh and Rendl \cite{PRendl09} have the same optimal value, but the optimal solution may be different except that the rank of the optimal solution is one. Because it is well-known that the rank constrained matrix optimization problems are computationally intractable and difficult in general \cite{BussFS99}, the rank one constraints are usually dropped in both the CPP and its related DNN relaxations of QAP. However, by use of the strategy of the difference of two convex functions (DC), the rank constraint can be replaced by the difference of the nuclear norm function and Ky-Fan $k$-norm function. Based on this simple observation, a penalty approach are proposed by \cite{Yan10} for calibrating rank constrained correlation matrix problems, which usually performances very well in many applications (see also \cite{LQi11}). In fact, based on the DC reformulations of the rank constraints, we shall reformulate the original QAP as a DC programming \cite{LTH12, LTao18} and employ the DC algorithm (DCA) to solve the non-convex QAP relaxation problems.

In this paper, we will propose a new rank constrained DNN model and show that it is equivalent with the original QAP (in the sense of both optimal values and optimal solutions). Also, we shall show the same techniques can be applied by other important non-convex problems such as the standard quadratic programming and the minimum-cut graph tri-partitioning problem. Although the equivalent rank constrained DNN model is still numerically intractable, we will propose a semi-proximal DC algorithm (DCA) framework for finding a feasible stationary point. Furthermore, for the large-scaled DCA inner subproblems, we will apply an efficient majorized semismooth Newton-CG augmented Lagrangian method based on  the software package SDPNAL+ \cite{STYZ19}. Finally, numerical experiments on the QAPLIB \cite{HAQAPLIB}  and  `dre' instances \cite{DreznerHT05} demonstrate the proposed approach usually performs well.

Below are some common notations to be used in this paper. We use $\mathcal{S}^{q}$ to denote the linear subspace of all $q \times q$ real symmetric matrices. Let $\mathcal{N}^q\subseteq\mathcal{S}^{q}$ be the subset of all $q \times q$ nonnegative symmetric matrices in $\mathcal{S}^{q}$. Denote  $\mathcal{S}^{q}_+/\mathcal{S}^{q}_-$ ($\mathcal{S}^{q}_{++}/\mathcal{S}^{q}_{--}$) the positive/negative semidefinite (definite) matrix cone in $\mathcal{S}^{q}$.  Moreover, let ${\cal C}^{q}$ be the set of copositive matrices in ${\cal S}^q$ and $({\cal C}^{q})^*$ be the dual cone of ${\cal C}^{q}$, i.e., the set of all completely positive matrices in ${\cal S}^q$. For a given matrix $Z\in {\cal S}^{q^2}$ with $q\ge 1$, we also use the following block notation for simplicity:
\begin{equation*}
	Z=\left[ \begin{array}{ccc}
		Z^{11} & \cdots & Z^{1q} \\ [3pt]
		\vdots & \ddots & \vdots \\ [3pt]
		Z^{q1} & \cdots & Z^{qq}
	\end{array} \right]
\end{equation*}
with $Z^{ij}\in\mathcal{R}^{q\times q}$ for each $i,j\in \{1,\ldots,q\}$. Let  $e_i$ be the $i$-th standard unit vector. We denote the vector and square matrix
of all ones by ${\bf 1}_q$ and $E_q$ respectively, and denote the identity matrix by $I_q$. We will omit the superscript $q$ if the dimension is clear. For a given $Z\in{\cal S}^q$, we use $\lambda_1(Z)\ge\ldots\ge \lambda_q(Z)$ to denote the eigenvalues of $Z$ (all real and counting multiplicity) arranging in non-increasing order.
We use ``${\rm vec}(\cdot)$" to denote the vectorization of matrices and use ``${\rm mat}(\cdot)$" to denote its inverse operator, i.e., the corresponding matricization of vectors. If $z\in {\cal R}^q$, then ${\rm Diag}(z)$ is a $q\times q$ diagonal matrix with $z$ on the main diagonal. Finally, we use  ``$\otimes$" to denote the Kronecker product between matrices.

\section{The rank constrained DNN reformulation of the QAP}\label{section:QAP-conic}

It is well-known that each permutation $\pi\in {\cal P}^n$ can be represented by a $n\times n$ permutation matrix $X$, i.e., a square binary matrix which has exactly one entry of $1$ in each row and each column and zeros elsewhere. Therefore, the QAP \eqref{eq:QAP-original} can be reformulate as the following trace form:
\begin{equation}\label{eq:QAP-trace}
	\min\left\{ \langle X, AXB+C\rangle \mid X\in \Pi^{n\times n} \right\},
\end{equation}
where $\langle \cdot,\cdot\rangle$ stands for the standard trace inner product of matrices, i.e., $\langle Y,Z\rangle={\rm tr}(X^TY)$ for $X,Y\in \mathcal{R}^{m\times n}$, and $\Pi^{n\times n}$ is the set of all $n\times n$ permutation matrices. It is clear that $\Pi^{n\times n}$ be characterized by the interaction of the set of orthogonal matrices and the set of nonnegative matrices, i.e.,
$$
\Pi^{n\times n}=\left\{X\in \mathcal{R}^{n\times n}\mid X^TX=I,\ X\ge 0 \right\}.
$$
Without loss of generality, we may assume that the data matrices $A,B,C$ in \eqref{eq:QAP-original} are nonnegative, i.e.,  $A,B,C\in \mathcal{N}$. Inspired by \cite{AWolkowicz00}, Povh and Rendl \cite{PRendl09} suggested to consider the following convex completely positive conic relaxation  of the QAP \eqref{eq:QAP-trace}:
\begin{equation}\label{eq:QAP_CP}
	\begin{array}{cl}
\displaystyle\min &  \langle B\otimes A+{\rm Diag}(c), Y \rangle \\[3pt]
\text{s.t.}  & \displaystyle\sum_{i=1}^nY^{ii}=I, \quad  \langle I, Y^{ij} \rangle = \delta_{ij}, \quad i,j\in \{1,\ldots,n \}, \\ [3pt]
& \langle E,Y \rangle = n^2, \quad Y\in({\cal C}^{n^2})^*,
\end{array}
\end{equation}
where $c={\rm vec}(C)$ and $\delta_{ij}=1$ if $i=j$ and $\delta_{ij}=0$ otherwise for $i,j\in\{1,\ldots,n\}$. It is clear that for any $n\times n$ permutation matrix $X\in \Pi^{n\times n}$,
\begin{equation}\label{eq:rank1-Y}
Y={\rm vec}(X){\rm vec}(X)^T, \quad X\in \Pi^{n\times n}.
\end{equation} is a feasible solution of \eqref{eq:QAP_CP}. Furthermore, Povh and Rendl \cite{PRendl09}  shown that the optimal value of \eqref{eq:QAP_CP} is actually equal the optimal value of QAP \eqref{eq:QAP-trace}. Unfortunately, the completely positive cone constrain $Y\in({\cal C}^{n^2})^*$ is computational intractable. A useful strategy to handle this is to approximate the cone $({\cal C}^{n^2})^*$ from the outside, e.g., the cone of symmetric doublely nonnegative matrices ${\cal S}^{n^2}_+\bigcap {\cal N}^{n^2}$. Thus, we obtain the following relaxation of the QAP \eqref{eq:QAP-trace}:
\begin{equation}\label{eq:QAP_DN}
\begin{array}{cl}
\displaystyle\min &  \langle B\otimes A+{\rm Diag}(c), Y \rangle \\[3pt]
\text{s.t.}  & \displaystyle\sum_{i=1}^nY^{ii}=I, \quad  \langle I, Y^{ij} \rangle = \delta_{ij}, \quad i,j\in \{1,\ldots,n \}, \\ [3pt]
& \langle E,Y \rangle = n^2, \quad Y\in {\cal S}^{n^2}_+\bigcap {\cal N}^{n^2}.
\end{array}
\end{equation}
Clearly, the optimal value of  problem \eqref{eq:QAP_DN} only provides a lower bound of the QAP \eqref{eq:QAP-trace}. In general, the relaxation \eqref{eq:QAP_DN} for the QAP is not tight.

On the other hand, from the equation \eqref{eq:rank1-Y}, we may add the rank constraint ${\rm rank}(Y)\le  1$ to \eqref{eq:QAP_DN} and obtain the following rank constrained doubly nonnegative (DNN) problem:
\begin{equation}\label{eq:QAP_SDP_rank}
\begin{array}{cl}
\displaystyle\min &  \langle B\otimes A+{\rm Diag}(c), Y \rangle \\[3pt]
\text{s.t.}  & \displaystyle\sum_{i=1}^nY^{ii}=I, \quad \langle I, Y^{ij} \rangle = \delta_{ij}, \quad i,j\in \{1,\ldots,n \}, \\ [3pt]
& \langle E,Y \rangle = n^2, \quad Y\in{\cal S}^{n^2}_+\bigcap {\cal N}^{n^2},  \quad {\rm rank}(Y)\le 1.
\end{array}
\end{equation}
The resulting problem \eqref{eq:QAP_SDP_rank} is non-convex. In fact,  we shall show that \eqref{eq:QAP_SDP_rank} is an exact reformulation of the original QAP \eqref{eq:QAP-trace}. To this end, we need the following simple observation on the rank one completely positive matrices.

\begin{lemma}\label{lem:rankone-CP}
	Let $q\ge 1$ be a given positive integer. Suppose that  $Y\in {\cal S}^{q}$ and  ${\rm rank}(Y)\le 1$. Then, the following statements are equivalent:
	\begin{itemize}
		\item[(i)]  $Y\in ({\cal C}^{q})^*$;
		
		\item[(ii)] $Y\in {\cal S}^{q}_+\bigcap {\cal N}^{q}$;
		
		\item[(iii)] there exists $x\in \mathcal{R}^{q}_+$ such that $Y=xx^T$.
	\end{itemize}
\end{lemma}
\begin{proof}
	Since ``(i) $\Longrightarrow$ (ii)" and ``(iii) $\Longrightarrow$ (i)" are obvious, we only need to show ``(ii) $\Longrightarrow$ (iii)", i.e., if $Y\in {\cal S}^{q}_+\bigcap {\cal N}^q$, then there exists $x\in \mathcal{R}^{q}_+$ such that $Y=xx^T$. Without loss of generality, we may assume ${\rm rank}(Y) = 1$, since otherwise the result holds trivially. It follows from $Y\in {\cal S}^{q}_+$ and ${\rm rank}(Y)=1$ that there exists $u\in \mathcal{R}^{q}$ such that $Y=\lambda uu^T$. Since $Y\ge 0$,  we have $Y_{ij}=u_iu_j\ge 0$ for each $i,j\in\{1,\ldots,q\}$. Thus, we can choose $x=\sqrt{\lambda} u \in \mathcal{R}^q_+$ such that $Y=xx^T$. \qed
\end{proof}

It is clear that the objective functions of \eqref{eq:QAP-trace} and \eqref{eq:QAP_SDP_rank} coincide. The equivalence between \eqref{eq:QAP-trace} and \eqref{eq:QAP_SDP_rank} then follows if we show the feasible sets of these two problems are the same. By employing the similar argument as that of \cite[Theorem 3]{PRendl09}, we have the following result on the equivalence of the feasible sets of \eqref{eq:QAP_SDP_rank} and \eqref{eq:QAP-trace}.

%
%

\begin{proposition}\label{prop:cp-rank-f-eq}
The matrix $Y\in {\cal S}^{n^2}_+$ is a feasible solution of \eqref{eq:QAP_SDP_rank} if and only if there exists a unique $X\in \Pi^{n\times n}$ such that $Y= {\rm vec}(X){\rm vec}(X)^T$. Moreover, since $\|{\rm vec}(X)\|^2$ is the only nonzero eigenvalue of $Y$, the vector ${\rm vec}(X)/ \|{\rm vec}(X)\|$ is the unit nonnegative eigenvector of $Y$.
\end{proposition}
\begin{proof}
	It is easy to see that if $X\in\Pi^{n\times n}$ then $Y={\rm vec}(X){\rm vec}(X)^T$ belongs the feasible set of \eqref{eq:QAP_SDP_rank}. Thus, we only need to show the converse direction holds. Suppose that $Y$ is a feasible set of  \eqref{eq:QAP_SDP_rank}. We know that ${\rm rank}(Y)=1$, since $Y\neq 0$. It then follows from Lemma \ref{lem:rankone-CP} that there exists $y \in \mathcal{R}^{n^2}_+$ such that $Y=yy^T$. Denote $X={\rm mat}(x)\in \mathcal{R}^{n\times n}$. Then, by employing the similar argument as that of \cite[Theorem 3]{PRendl09}, we are able to show that $X\in \Pi^{n\times n}$. Furthermore, it is easy to verify  that for any $X,X'\in\Pi^{n\times n}$, if $X\neq X'$, then $Y\neq Y'$ with $Y={\rm vec}(X){\rm vec}(X)^T$ and $Y'={\rm vec}(X'){\rm vec}(X')^T$.
	
	Let the nonzero unit vector $v \in \mathcal{R}^{n^2}$ with $v ={\rm vec}(X)/\|{\rm vec}(X)\| $, Obviously, $v \in \mathcal{R}^{n^2}_+$. From the definition of the characteristic polynomial for matrices, we know that
	$$ Y v = {\rm vec}(X){\rm vec}(X)^T \cdot  \frac{{\rm vec}(X) }{\|{\rm vec}(X)\|}= \|{\rm vec}(X)\|  {\rm vec}(X) = \|{\rm vec}(X)\|^2 v,$$
	that is , $\|{\rm vec}(X)\|^2$ and $v$ is the eigenvalue and eigenvector of $Y$ respectively. 
The proof is completed. \qed
\end{proof}

\vskip 3 true pt

\begin{remark}
It  follows from  Proposition \ref{prop:cp-rank-f-eq} that if $Y\in\mathcal{S}^{n^2}$ is a feasible solution of \eqref{eq:QAP_SDP_rank}, then we can find the permutation matrix $X\in\Pi^{n\times n}$ by setting $X={\rm mat}(x)$ with $x = v \cdot \|{\rm vec}(X)\|$ easily, where $v$ is the unit corresponding eigenvector of $Y$ with respect to $n$. 
\end{remark}

\vskip 3 true pt

The following result on the equivalence between the rank constrained DNN problem \eqref{eq:QAP_SDP_rank} and the QAP \eqref{eq:QAP-trace} follows from Proposition \ref{prop:cp-rank-f-eq} immediately.

\vskip 4 true pt

\begin{theorem}\label{thm:QAP-SDP}
	The rank constrained DNN problem \eqref{eq:QAP_SDP_rank} is equivalent to the QAP \eqref{eq:QAP-trace}.
\end{theorem}

\vskip 3 true pt
Clearly, the non-convex rank constrained DNN representation  \eqref{eq:QAP_SDP_rank} is at least as hard as the original QAP, which means that finding a global solution of \eqref{eq:QAP_SDP_rank} is computational intractable. However, it is still possible to design some efficient algorithms, e.g., the DCA (see Section \ref{section:DCA}), to find a good feasible point of \eqref{eq:QAP_SDP_rank} and obtain a good feasible solution of the original QAP.

\section{Extensions}
In this section, we shall demonstrate that the results obtained in Section \ref{section:QAP-conic} can be applied to other important non-convex problems, which have the similar rank constrained DNN representations.

{\bf Standard quadratic programming}. The standard quadratic problem (StQP) consists of finding an optimal of a quadratic form over the standard simplex, i.e.,
\begin{equation}\label{eq:StQP}
\min\Big\{ \left\langle x, Qx\right\rangle \mid \sum_{i=1}^nx_i=1, \  x\ge 0 \Big\},
\end{equation}
where $Q$ is an arbitrary $n\times n$ symmetric matrix. The StQP \eqref{eq:StQP} includes many important combinatorial optimization problems as special cases, e.g., the maximum clique problem \cite{MStraus65}.  It is clear that the StQP \eqref{eq:StQP} can be rewritten as the following matrix form:
  \begin{equation*}
  \begin{array}{cl}
  \displaystyle\min & \langle Q, Y \rangle \\[3pt]
  \text{s.t.}  &  \langle E, Y \rangle = 1, \quad Y=xx^T,\quad  x\ge 0.
  \end{array}
  \end{equation*}
  Thus, by employing Lemma \ref{lem:rankone-CP}, we obtain the following result on the rank constrained DNN representation of  the StQP \eqref{eq:StQP}, immediately.
\begin{theorem}
	The standard quadratic problem \eqref{eq:StQP} is equivalent to the following rank constrained DNN problem:
	 \begin{equation}\label{eq:StQP_rank}
	\begin{array}{cl}
	\displaystyle\min & \langle Q, Y \rangle \\[3pt]
	\text{s.t.}  &  \langle E, Y \rangle = 1, \quad Y\in{\cal S}^n_+\bigcap {\cal N}^n,\quad  {\rm rank}(Y)\le 1.
	\end{array}
	\end{equation}
\end{theorem}

{\bf The minimum-cut graph tri-partitioning problem}. The minimum-cut graph tri-partitioning problem \cite{PRendl07} is to find a tri-partitioning the vertices of a graph into sets $S_1$, $S_2$ and $S_3$ of specified
cardinalities, such that the total weight of edges between $S_1$ and $S_2$ is minimal.

Let $G = (V, E)$ be an undirected graph on $n$ vertices, given by its (weighted) symmetric nonnegative adjacency matrix $A \in \mathcal{N}^{n}$,
the minimum-cut graph tri-partitioning problem \cite{PRendl07}  can be described as:
for given integers $m_1$, $m_2$ and $m_3$ summing to $n$,
find subsets $S_1$, $S_2$ and $S_3$ of $V(G)$ with cardinalities $m_1$, $m_2$ and $m_3$, respectively,
such that the total weight of edges between $S_1$ and $S_2$ is minimal. By presenting partitions $S_1$, $S_2$ and $S_3$ by $n\times 3$ matrices $X$, the minimum-cut graph tri-partitioning problem can be written as follows
\begin{equation}\label{eq:mat_form_mincut}
\begin{array}{cl}
\displaystyle\min &  \displaystyle\frac{1}{2}\langle X, AXB\rangle \\[3pt]
\text{s.t.}  & X^T X = M, \quad X {\bf 1}_3 = {\bf 1}_n, \\ [3pt]
 & X\ge 0,
\end{array}
\end{equation}
where $M:={\rm Diag}(m_1,m_2,m_3)$ and $B=\left[
\begin{array}{ccc}
0 & 1 & 0 \\
1 & 0 & 0 \\
0 & 0 & 0
\end{array}
\right]$, the vector of all ones is $ {\bf 1}_k \in \mathcal{R}^k$.
By introducing $Y=xx^T$ with $x={\rm vec}(X)$,
Povh and Rendl \cite{PRendl07} reformulate the minimum-cut graph tri-partitioning problem \eqref{eq:mat_form_mincut} as follows:
\begin{equation}\label{eq:mincut}
\begin{array}{cl}
\displaystyle\min &  \displaystyle\frac{1}{2}\langle B\otimes A, Y \rangle \\[3pt]
\text{s.t.}  & \langle L^{ij}\otimes I, Y \rangle = m_i\delta_{ij}, \quad 1\leq i\leq j\leq 3, \\ [3pt]
& \langle E_3\otimes J^{ii},Y \rangle = 1, \quad 1\leq i\leq n, \\ [3pt]
& \langle V_i\otimes W_{j}^T,Y \rangle = m_i, \quad 1\leq i\leq 3, 1\leq j\leq n \\ [3pt]
& \langle L^{ij}\otimes E_n,Y \rangle = m_i m_j, \quad 1\leq i\leq j\leq 3, \\ [3pt]
& Y = xx^T, x\in\mathcal{R}^{3n}_+,
\end{array}
\end{equation}
where $V_i=e_i {\bf 1}_3^T\in\mathcal{R}^{3\times 3}$ for $i=1,2,3$, $W_{j}=e_j {\bf 1}_n^T\in \mathcal{R}^{n\times n}$ for $j=1,\ldots,n$, $J^{ij}=e_i e_j^T\in\mathcal{R}^{n\times n}$ and $L^{ij}=\frac{1}{2}(e_ie_j^T+e_je_i^T)\in\mathcal{R}^{3\times 3}$ for $i,j=1,2,3$. Again, similar with Section \ref{section:QAP-conic}, by employing Lemma \ref{lem:rankone-CP}, we are able to obtain the following rank constrained DNN representation of  the minimum-cut graph tri-partitioning problem \eqref{eq:mat_form_mincut}.
\begin{theorem}
	The minimum-cut graph tri-partitioning problem \eqref{eq:mat_form_mincut} is equivalent to the following rank constrained DNN problem:
\begin{equation}\label{eq:mincut-DNN}
\begin{array}{cl}
\displaystyle\min &  \displaystyle\frac{1}{2}\langle B\otimes A, Y \rangle \\[3pt]
\text{s.t.}  & \langle L^{ij}\otimes I, Y \rangle = m_i\delta_{ij}, \quad 1\leq i\leq j\leq 3, \\ [3pt]
& \langle E_3\otimes J^{ii},Y \rangle = 1, \quad 1\leq i\leq n, \\ [3pt]
& \langle V_i\otimes W_{j}^T,Y \rangle = m_i, \quad 1\leq i\leq 3, 1\leq j\leq n \\ [3pt]
& \langle L^{ij}\otimes E_n,Y \rangle = m_i m_j, \quad 1\leq i\leq j\leq 3, \\ [3pt]
& Y\in{\cal S}^{3n}_+\bigcap {\cal N}^{3n},\quad  {\rm rank}(Y)\le 1.
\end{array}
\end{equation}
\end{theorem}

\section{The DCA for the rank constrained DNN problem}\label{section:DCA}

In this section, we shall propose a DCA based algorithm for the rank constrained DNN relaxations established in the previous section. For simplicity in notation, all proposed rank constrained DNN representations \eqref{eq:QAP_SDP_rank}, \eqref{eq:StQP_rank} and \eqref{eq:mincut-DNN} can be cast in the following
abstract form:
\begin{equation}\label{eq:QAP_SDP_rank_simp}
\begin{array}{cl}
\displaystyle\min &  f(Y):=\langle \overline{C}, Y \rangle \\[3pt]
\text{s.t.}  & Y\in \Omega\bigcap {\cal R},
\end{array}
\end{equation}
where the subsets $\Omega,{\cal R}\subseteq {\cal S}^q$ are defined  by
\begin{equation}\label{eq:def-Omega}
	\Omega:=\left\{Y\in{\cal S}^q_+\bigcap {\cal N}^q \mid {\cal A}(Y)=b\right\}
\end{equation}
and
\begin{equation}\label{eq:def-R-set}
{\cal R}:=\left\{ Y\in{\cal S}^{q}\mid {\rm rank}(Y)\le 1 \right\},
\end{equation}
$\overline{C}\in {\cal S}^q$, ${\cal A}:{\cal S}^q\to {\cal R}^m$ is a given linear operator, and $b\in {\cal R}^{m}$ is a given data.

It is worth to note that for the rank constrained DNN relaxations proposed in Section \ref{section:QAP-conic}, the subsets $\Omega$ with respect to \eqref{eq:QAP_SDP_rank}, \eqref{eq:StQP_rank} and \eqref{eq:mincut-DNN} are satisfy the following assumption.

\begin{assumption}\label{ass:bounded}
	The subset $\Omega\subseteq{\cal S}^q$ defined by \eqref{eq:def-Omega} is nonempty and bounded.
\end{assumption}

Let $\rho >0$ be a given penalty parameter. The rank constrained DNN problem \eqref{eq:QAP_SDP_rank_simp} is closed related to the following rank penalized problem:
\begin{equation}\label{eq:QAP_SDP_rank_p}
\begin{array}{cl}
\displaystyle\min  &   f(Y) +\rho\, {\rm rank}(Y) \\[3pt]
\text{s.t.}  & Y\in \Omega.
\end{array}
\end{equation}
In fact, we shall verify  that under Assumption \ref{ass:bounded}, the rank penalized problem \eqref{eq:QAP_SDP_rank_p} is an exact penalty version of the rank constrained DNN problem \eqref{eq:QAP_SDP_rank_simp} in the sense that there exists a constant $\overline{\rho}>0$ such that the global optimal solution of \eqref{eq:QAP_SDP_rank_p} associated to any $\rho\ge\overline{\rho}$ coincides with that of \eqref{eq:QAP_SDP_rank_simp}.

\begin{theorem}\label{thm:exact-p-rank}
	Suppose Assumption \ref{ass:bounded} holds. There exists a constant $\overline{\rho}>0$ such that for any $\rho\ge \overline{\rho}$, the global optimal solution set of \eqref{eq:QAP_SDP_rank_p} associated to any $\rho >\overline{\rho}$ coincides with the global optimal solution set  of \eqref{eq:QAP_SDP_rank_simp}.
\end{theorem}
\begin{proof}
	Let $Y^*$ be a global optimal solution of \eqref{eq:QAP_SDP_rank_simp}. Since $\Omega$ is assumed nonempty and compact, we may assume that $\widetilde{Y}\in \Omega$ is an optimal solution of the convex problem $\min \left\{ f(Y) \mid Y\in \Omega  \right\}$.   It is clear that $f(Y^*)\ge f(\widetilde{Y})$. Let $\overline{\rho}>f(Y^*)- f(\widetilde{Y})\ge 0$ be fixed. Suppose that $\rho\ge \overline{\rho}$. Let $Y_{\rho}$ be a global optimal solution of \eqref{eq:QAP_SDP_rank_p} chosen arbitrarily with respect to $\rho$. We have
	\begin{equation}\label{eq:rank-ineq-2}
	f(Y_{\rho})+\rho\, {\rm rank}(Y_{\rho})\le f(Y^*)+\rho\, {\rm rank}(Y^*)\le f(Y^*)+\rho.
	\end{equation}
	By noting that ${\rm rank}(Y_{\rho})\ge 1$ (since $Y_{\rho}\neq 0$), we obtain from \eqref{eq:rank-ineq-2} that
	\begin{equation}\label{eq:f-ineq-exact-rank}
		f(Y_{\rho})\le  f(Y^*).
	\end{equation}
	Since $Y_{\rho}\in \Omega$, we have $f(\widetilde{Y})\le f(Y_{\rho})$. Thus, we have
	\begin{equation}\label{eq:rank-ineq-1}
	\rho({\rm rank}(Y_{\rho})-1)\le f(Y^*)-f(\widetilde{Y}).
	\end{equation}
	We claim that ${\rm rank}(Y_{\rho})\le 1$. In fact, if ${\rm rank}(Y_{\rho})\ge 2$, then it follows from  \eqref{eq:rank-ineq-1} that
	\[
	\rho\le f(Y^*)-f(\widetilde{Y}),
	\]
	which contradicts with the fact that $\rho\ge \overline{\rho}>f(Y^*)- f(\widetilde{Y})$. Thus, we know that $Y_{\rho}\in \Omega\bigcap {\cal R}$, i.e., $Y_{\rho}$ is indeed a feasible solution of \eqref{eq:QAP_SDP_rank_simp}. Therefore, we have  $f(Y_{\rho})\ge f(Y^*)$ since $Y^*$ is a global solution  of \eqref{eq:QAP_SDP_rank_simp}. This, together with \eqref{eq:f-ineq-exact-rank}, implies that $f(Y_{\rho})= f(Y^*)$, which implies that $Y_{\rho}$ is  a global solution of \eqref{eq:QAP_SDP_rank_simp}. On the other hand, by noting that $Y_{\rho}\neq 0$ and ${\rm rank}(Y_{\rho})\le 1$, we conclude that ${\rm rank}(Y_{\rho})= 1$, which implies that
	\[
	f(Y_{\rho})+\rho\, {\rm rank}(Y_{\rho})=f(Y_{\rho})+\rho \ge  f(Y^*)+\rho\, {\rm rank}(Y^*)
	\]
	It then follows from  \eqref{eq:rank-ineq-2} that $f(Y_{\rho})+\rho\, {\rm rank}(Y_{\rho})=f(Y^*)+\rho\, {\rm rank}(Y^*)$. Thus, we know that $Y^*$ is also a global solution of \eqref{eq:QAP_SDP_rank_p}. Since $Y^*$ and $Y_{\rho}$ are chosen arbitrarily, we know that the global solution sets of \eqref{eq:QAP_SDP_rank_simp} and \eqref{eq:QAP_SDP_rank_p} coincide. \qed
\end{proof}

Consider the following penalized problem:
\begin{equation}\label{eq:QAP_SDP_p}
\begin{array}{cl}
\displaystyle\min  &   f_{\rho}(Y):=\langle \overline{C}, Y \rangle +\rho (\|Y\|_*-\|Y\|_2) \\[3pt]
\text{s.t.}  & Y\in \Omega.
\end{array}
\end{equation}

Let $\mathcal{X}$ and $\mathcal{Z}$ be two finite dimensional Euclidean space. Recall a set-valued mapping $\Psi:{\cal X}\rightrightarrows \mathcal{Z}$ is called calm at $\bar{x}$ for $\bar{z}\in\Psi(\bar{x})$ if there exist a constant $\alpha>0$ and neighborhood ${\cal U}\subseteq \mathcal{X}$ of $\bar{x}$ and neighborhood ${\cal V}\subseteq \mathcal{Z}$ of $\bar{z}$ such that
\begin{equation*}
\Psi(x)\cap {\cal V}\subseteq \Psi(\bar{x})+\alpha\|x-\bar{x}\|\mathbb{B}_{\cal X}\quad \forall\, x\in \mathcal{U},
\end{equation*}
where $\mathbb{B}_{\cal X}$ is the unit ball in $\mathcal{X}$.

\begin{proposition}\label{prop:exact-p}
Suppose that the set-valued mapping $\Gamma:\Re\rightrightarrows \mathcal{S}^q $ defined by
	\[
	\Gamma(w):=\left\{ Y\in \mathcal{S}^q \mid Y\in \Omega,\ \|Y\|_*-\|Y\|_2=w \right\},\quad w\in\Re,
	\] is calm at $0$ for each $Y\in \Gamma(0)$. Then, there exists a constant $\bar{\rho}>0$ such that for any $\rho>\bar{\rho}$, $Y^*$ is an optimal of \eqref{eq:QAP_SDP_rank_simp} if and only if $Y^*$ is an optimal of the penalized problem \eqref{eq:QAP_SDP_p}.
\end{proposition}
\begin{proof}
	First, we shall show that there exists $\bar{\rho}>0$  if $Y^*$ is an optimal of \eqref{eq:QAP_SDP_rank_simp}, then it is also an optimal solution of the  penalized problem \eqref{eq:QAP_SDP_p} for $\rho>\bar{\rho}$. By \cite[Theorem 2.1]{BPan16}, we know from the calmness of $\Gamma$ that there exists $\tau>0$ such that ${\rm dist}(\overline{Y},\Omega\bigcap {\cal R})\le \tau {\rm dist}(\overline{Y}, {\cal R})=\tau(\|\overline{Y}\|_*-\|\overline{Y}\|_2)$. Let $L:=\|\overline{C}\|>0$. Suppose that $\rho>\bar{\rho}:=\max\{L\tau,L \}$ be arbitrarily given.  Suppose  there exists $\overline{Y}\in\Omega$ and $\varepsilon>0$ such that
	\[
	 \langle \overline{C}, \overline{Y} \rangle + \rho (\|\overline{Y}\|_*-\|\overline{Y}\|_2) < \langle \overline{C}, Y^* \rangle - \rho \varepsilon.
	\]
	Let $\widehat{Z}\in \Omega\bigcap {\cal R}$ be such that
	\[
	\|\widehat{Z}-\overline{Y}\|\le {\rm dist}(\overline{Y},\Omega\bigcap {\cal R})+\varepsilon.
	\]
	Since ${\rm dist}(\overline{Y},\Omega\bigcap {\cal R})\le \tau {\rm dist}(\overline{Y}, {\cal R})=\tau(\|\overline{Y}\|_*-\|\overline{Y}\|_2)$.  we have
	\[
	\|\widehat{Z}-\overline{Y}\|\le\tau( \|\overline{Y}\|_*-\|\overline{Y}\|_2)+\varepsilon.
	\]
	 Then, 
	\begin{eqnarray*}
	\langle \overline{C}, \widehat{Z} \rangle &\le& \langle \overline{C}, \overline{Y} \rangle + L\|\widehat{Z}-\overline{Y}\| \le \langle \overline{C}, \overline{Y} \rangle + L(\tau( \|\overline{Y}\|_*-\|\overline{Y}\|_2)+\varepsilon) \\
	&\le&  \langle \overline{C}, \overline{Y} \rangle + \rho( \|\overline{Y}\|_*-\|\overline{Y}\|_2+\varepsilon) < \langle \overline{C}, Y^* \rangle.
	\end{eqnarray*}
This contradicts with the fact that $Y^*$ is an optimal of \eqref{eq:QAP_SDP_rank_simp}.

For the converse direction, it is sufficient to show that if $Y^*$ is an optimal of the penalized problem \eqref{eq:QAP_SDP_p}, then $Y^*\in \Omega\bigcap {\cal R}$, i.e.,  $Y^*$ is a feasible solution of \eqref{eq:QAP_SDP_rank_simp}. In fact, if $\widetilde{Y}\in \Omega\bigcap {\cal R}$ is an optimal of  \eqref{eq:QAP_SDP_rank_simp}, then since $Y^*$ is an optimal of the problem \eqref{eq:QAP_SDP_p}, we know from the first part that
\[
\langle \overline{C}, Y^* \rangle + \rho (\|Y^*\|_*-\|Y^*\|_2)=\langle \overline{C}, \widetilde{Y} \rangle
\]
and
\[
\langle \overline{C}, Y^* \rangle + \frac{1}{2}(\rho+\bar{\rho}) (\|Y^*\|_*-\|Y^*\|_2)\ge \langle \overline{C}, \widetilde{Y} \rangle,
\]
which implies that
\[
\frac{1}{2}(\bar{\rho}-\rho) (\|Y^*\|_*-\|Y^*\|_2)\ge 0.
\]
Since $\rho>\bar{\rho}$ and $\|Y^*\|_*-\|Y^*\|_2\ge 0$, we know that $\|Y^*\|_*-\|Y^*\|_2=0$, i.e., ${\rm rank}(Y^*)\le 1$. Thus, we have $Y^*\in \Omega\bigcap{\cal R}$. This completes the proof. \qed
\end{proof}

The objective function of \eqref{eq:QAP_SDP_p} can be rewritten as
\begin{eqnarray*}
f_{\rho}(Y)=\langle \overline{C},Y \rangle + \rho \|Y\|_* - \rho p(Y), \quad Y\in {\cal S}^q,
\end{eqnarray*}
where $p(Y):=\|Y\|_2$. Therefore, the non-convex objective function of the penalized problem  \eqref{eq:QAP_SDP_p} is a DC (difference of convex) function. Thus, we introduce a DC based algorithm to solve \eqref{eq:QAP_SDP_p}, which has the following template:
\begin{algorithm}\caption{[Proximal DC Algorithm (ProxDCA)]}\label{alg:DCA}
\begin{algorithmic}[1]
\STATE  Let $Y^0\in\Omega$ be an initial point and $\sigma>0$. Set $k=0$.
 \STATE Choose $W^{k}\in \partial\,p(Y^{k})$. Compute
     \begin{equation}\label{eq:prox_mod}
     Y^{k+1}= \argmin\left\{ \widehat{f}_{\rho,\sigma}(Y)  \mid Y\in\Omega \right\},
     \end{equation}
     where
  \begin{equation}\label{prox_sdp}
     \widehat{f}_{\rho,\sigma}(Y):=\langle \overline{C},Y \rangle + \rho \|Y\|_* - \rho( p(Y^{k})+\langle W^{k},Y-Y^{k} \rangle)+\displaystyle\frac{1}{2{\sigma}}\|Y-Y^{k}\|^2
 \end{equation}
      and the subset $\Omega\subseteq {\cal S}^q$ is defined by \eqref{eq:def-Omega}.

\STATE If $Y^{k+1}=Y^k$ stop; otherwise set $k=k+1$ and go to {\bf Step.2}.
\end{algorithmic}

\end{algorithm}

Under Assumption \ref{ass:bounded}, the strongly convex problem \eqref{eq:prox_mod} has a unique solution and can be solved efficiently by considering its dual problem, i.e.,
\begin{equation}\label{eq:prox_mod-dual}
\begin{array}{cl}
\displaystyle\max &  -\langle b,y \rangle -\displaystyle\frac{\sigma}{2}\left\|\overline{C}+\rho(I+W^k)+{\cal A}^*y+S+Z-\sigma Y^k\right\|^2 \\[3pt]
\text{s.t.}  & S\in{\cal S}^q_{-},\quad Z\in -{\cal N}^q.
\end{array}
\end{equation}
Moreover, if $(y^{k+1},S^{k+1},Z^{k+1})\in\mathcal{R}^m\times {\cal S}^q\times {\cal N}^q$ is an optimal solution of the above dual problem \eqref{eq:prox_mod-dual}, $Y^{k+1}$ can be found as follows
\begin{equation}
Y^{k+1}=Y^k-\sigma \left({\cal A}^*y^{k+1}+S^{k+1}+Z^{k+1}+\overline{C}+\rho(I-W^k)\right).
\end{equation}
It is clear that the dual problem \eqref{eq:prox_mod-dual} coincides with the inner problem \cite[(8)]{YangSunToh15} involved in the
augmented Lagrangian method of the dual problem of the semidefinite programming  with an additional polyhedral cone
constraint (SDP+) introduced by  \cite{YangSunToh15}. Therefore, we can employ the majorized semismooth Newton-CG method \cite[{\tt Algorithm MSNCG}]{YangSunToh15} to solve \eqref{eq:prox_mod}, directly. Furthermore, in order for the dual problem \eqref{eq:prox_mod-dual} to have a bounded solution set, we introduce the following  general Slater condition for the constraint set $\Omega$ defined in \eqref{eq:def-Omega}.

\begin{assumption}\label{ass:Slater}
	There exists  $\widetilde{Y}\in {\cal S}^q$ such that
	\[
	{\cal A}({\cal T}_{{\cal N}^q}(\widetilde{Y}))={\cal R}^m \quad {\rm and} \quad \widetilde{Y}\in {\cal S}^q_{++}\cap {\rm int}\,({\cal N}^q),
	\]
	where ${\rm int}\,({\cal N}^q)$ and ${\cal T}_{{\cal N}^q}(\widetilde{Y})$ denote the interior of ${\cal N}^q$ and the tangent cone of ${\cal N}^q$ at $\widetilde{Y}$, respectively.
\end{assumption}

Under Assumption \ref{ass:Slater}, the convergence of {\tt Algorithm MSNCG} is  established in \cite[Theorem 2.5]{YangSunToh15}. For simplicity, we omit details here.

Next, we shall study the convergence of the proposed DC based algorithm for the rank constrained DNN problem \eqref{eq:QAP_SDP_rank_simp}. A feasible point $Y\in \Omega$ is said to be a stationary point of the penalized problem \eqref{eq:QAP_SDP_p} if
\begin{eqnarray*}
	(\overline{C}+\rho I+\mathcal{N}_{\Omega}(Y))\bigcap(\rho\partial\,p(Y))\neq \emptyset,
\end{eqnarray*}
where $\mathcal{N}_{\Omega}(Y)$ is the normal cone of the convex set $\Omega$ at $Y$ in the sense of convex analysis (cf. e.g., \cite{Rockafellar70}).  We have the following results on the convergence of the proposed DC based algorithm (Algorithm \ref{alg:DCA}) for the rank constrained DNN problem \eqref{eq:QAP_SDP_rank_simp}. Note that the proof of the following proposition is similar with that of \cite[Theorem 3.4]{GSun10}. However, we include the proof here for completion.

\begin{proposition}\label{prop:convegence-DCA}
Suppose that Assumption \ref{ass:bounded} holds. Let $\rho>0$ be given. Let $\{Y^k\}$ be the sequence generated by Algorithm \ref{alg:DCA}. Then $\{f_{\rho}(Y^k)\}$ is a monotonically decreasing sequence.
If $Y^{k+1} = Y^k$ for some integer $k \geq 0$,
then $Y^{k+1}$ is a stationary point of the penalized problem \eqref{eq:QAP_SDP_p}.
Otherwise, the infinite sequence $\{f_{\rho}(Y^k)\}$ satisfies
\begin{eqnarray}\label{de_obj_fun}
\frac{1}{2{\sigma}}\|Y^{k+1} - Y^{k}\|^2 \leq f_{\rho}(Y^k) - f_{\rho}(Y^{k+1}),\quad k = 0,1,...
\end{eqnarray}
Moreover, any accumulation point of the bounded sequence $\{Y^k\}$ is a stationary point of
problem \eqref{eq:QAP_SDP_p}.
\end{proposition}
\begin{proof}
Since the function $p$ is convex and $W^k\in\partial\, p(Y^k)$, we know that
\begin{eqnarray*}
p(Y^{k+1}) &\geq& p(Y^k) + \langle W^k,Y^{k+1}-Y^k \rangle.
\end{eqnarray*}
Therefore, we have for each $k\ge 0$,
\begin{eqnarray}\label{decreasing_prop}
f_{\rho}(Y^{k+1})&=&\langle \overline{C},Y^{k+1} \rangle + \rho\|Y^{k+1}\|_* - \rho p(Y^{k+1}) \nn \\
 &\leq& \langle \overline{C},Y^{k+1} \rangle + \rho\|Y^{k+1}\|_* - \rho \big(p(Y^{k})+\langle W^k,Y^{k+1}-Y^k \rangle\big) \nn \\ 
 &&+\frac{1}{2{\sigma}}\|Y^{k+1}-Y^k\|^2  \leq \widehat{f}_{\rho,{\sigma}}(Y^k)=f_{\rho}(Y^k), \nn
\end{eqnarray}
where the last inequality due to $Y^k\in\Omega$ and $Y^{k+1}$ is the optimal solution of \eqref{eq:prox_mod}.
Thus, we know that the sequence $\{f_{\rho}(Y^k)\}$ is a monotonically decreasing sequence.

Assume that there exists some $k\ge 0$ such that $Y^{k+1}=Y^k$. We shall show that $Y^{k+1}$ is a stationary point of \eqref{eq:QAP_SDP_p}. Since $Y^{k+1}$ is the optimal solution of the strongly convex problem \eqref{eq:prox_mod},
we know that
\begin{eqnarray}\label{KKT_for_solvep}
0 \in \frac{1}{{\sigma}}(Y^{k+1}-Y^k)+\overline{C}-\rho W^k+\rho I+\mathcal{N}_{\Omega}(Y^{k+1}).
\end{eqnarray}
It then follows from $Y^{k+1}=Y^k$ that
\begin{eqnarray*}
\rho W^k \in \overline{C}+\rho I+\mathcal{N}_{\Omega}(Y^{k+1}),
\end{eqnarray*}
which implies that
\begin{eqnarray*}
(\overline{C}+\rho I+\mathcal{N}_{\Omega}(Y^{k+1}))\bigcap\rho\partial\,p(Y^{k+1})\neq \emptyset,
\end{eqnarray*}
i.e. $Y^{k+1}$ is a stationary point of \eqref{eq:QAP_SDP_p}.

Next, suppose that  for all $k\geq 0$, $Y^{k+1}\neq Y^{k}$.
It then follows from \eqref{KKT_for_solvep},
there exists $D^{k+1}\in\mathcal{N}_{\Omega}(Y^{k+1})$ such that
\begin{eqnarray}\label{eq:KKT_for_solvep}
0 = \frac{1}{{\sigma}}(Y^{k+1}-Y^{k})+\overline{C}-\rho(W^k-I)+D^{k+1}.
\end{eqnarray}
Thus, since $Y^{k}\in\Omega$ and $D^{k+1}\in\mathcal{N}_{\Omega}(Y^{k+1})$ for each $k\geq0$, by \eqref{decreasing_prop},
we have
\begin{eqnarray*}
&& f_{\rho}(Y^{k+1})-f_{\rho}(Y^{k}) \leq \widehat{f}_{\rho,{\sigma}}(Y^{k+1})-f_{\rho}(Y^{k}) \\
&=&\frac{1}{2{\sigma}}\|Y^{k+1}-Y^{k}\|^2+\langle \overline{C},Y^{k+1} \rangle - \rho(\|Y^k\|_2 + \langle W^k,Y^{k+1}-Y^k \rangle-\langle I,Y^{k+1} \rangle) \\
&& - (\langle \overline{C},Y^{k} \rangle - \rho(\|Y^k\|_2 -\langle I,Y^{k} \rangle) \\
&=&\frac{1}{2{\sigma}}\|Y^{k+1}-Y^{k}\|^2 + \langle \overline{C},Y^{k+1}-Y^{k} \rangle - \langle\rho(W^k-I),Y^{k+1}-Y^{k}\rangle \\
&=&\frac{1}{2{\sigma}}\|Y^{k+1}-Y^{k}\|^2 + \langle \overline{C}-\rho(W^k-I),Y^{k+1}-Y^{k} \rangle \\
&=&\frac{1}{2{\sigma}}\|Y^{k+1}-Y^{k}\|^2 + \langle -\frac{1}{{\sigma}}(Y^{k+1}-Y^{k})-D^{k+1},Y^{k+1}-Y^{k} \rangle \\
&=&-\frac{1}{2{\sigma}}\|Y^{k+1}-Y^{k}\|^2 - \langle D^{k+1},Y^{k+1}-Y^{k} \rangle \le 0,
\end{eqnarray*}
which implies that
\[
\frac{1}{2{\sigma}}\|Y^{k+1}-Y^{k}\|^2 \leq f_{\rho}(Y^{k})-f_{\rho}(Y^{k+1}).
\]
Thus, the infinite sequence $\{f_{\rho}(Y^k)\}$ satisfies the inequality \eqref{de_obj_fun}.

Moreover, suppose that $\overline{Y}$ is an accumulation point of $\{Y^k\}$.
Let $\{Y^{k_j}\}$ be a subsequence of $\{Y^k\}$ such that
$$
\lim_{j\rightarrow+\infty}Y^{k_j}=\overline{Y}.
$$
Then, by \eqref{de_obj_fun},  we obtain that
\begin{eqnarray*}
\lim_{i\to \infty}\frac{1}{2\sigma}\sum_{k=0}^i\|Y^{k+1}-Y^k\|^2 \le \liminf_{i\to \infty}(f_{\rho}(Y^0)-f_{\rho}(Y^{i+1}))\leq f_{\rho}(Y^0) < +\infty,
\end{eqnarray*}
which implies that $\displaystyle\lim_{k \to \infty}\|Y^{k+1}-Y^k\|=0$. Therefore, we obtain that
\[
\displaystyle\lim_{j\to\infty}Y^{k_j+1}=\lim_{j\to\infty}Y^{k_j}=\bar{Y} \quad {\rm and} \quad
\displaystyle\lim_{j\to\infty}(Y^{k_j+1}-Y^{k_j})=0.
\]
Furthermore, since $\{Y^{k_j}\}$ is bounded, it follows from \cite[Theorem 24.7]{Rockafellar70} that $\{W^{k_j}\}$ is also bounded.
By taking a subsequence if necessary, we may assume that there exists $\overline{W}\in\partial\, p(\overline{Y})$ such that $\displaystyle\lim_{j\to\infty}W^{k_j}=\overline{W}$.
Therefore, we obtain from \eqref{eq:KKT_for_solvep} that
\begin{eqnarray*}
\overline{D}:=\displaystyle\lim_{j\to\infty}D^{k_j+1}=\displaystyle\lim_{j\to\infty}-(\frac{1}{{\sigma}}(Y^{k_j+1}-Y^{k_j})+\overline{C}-\rho(W^{k_j}-I))=-\overline{C}-\rho I+\rho\overline{W}.
\end{eqnarray*}
Now in order to show that $\overline{Y}$ is a stationary point of problem \eqref{eq:QAP_SDP_p}, we only need to show that
$\overline{D}\in\mathcal{N}_{\Omega}(\overline{Y})$.
Suppose that $\overline{D}\notin\mathcal{N}_{\Omega}(\overline{Y})$, i.e., there exists $\widetilde{Y}\in\Omega$ such that $\langle\overline{D},\widetilde{Y}-\overline{Y}\rangle>0$.
Since for each $k_j$, $D^{k_j+1}\in\mathcal{N}_{\Omega}(Y^{k_j+1})$,
we have
\begin{eqnarray*}
\langle D^{k_j+1},\widetilde{Y}-Y^{k_j+1}\rangle \leq 0.
\end{eqnarray*}
It follows from the convergence of the two subsequences $\{D^{k_j+1}\}$ and $\{Y^{k_j+1}\}$,
thus
\begin{eqnarray*}
\langle \overline{D},\widetilde{Y}-\overline{Y}\rangle \leq 0.
\end{eqnarray*}
This is a contradiction. The proof is completed. \qed
\end{proof}

In order to show the infinity sequence $\{Y^k\}$ generated by the proposed Algorithm \ref{alg:DCA} actually converge, we recall the following definition of the Kurdyka-Łojaziewicz (KL) property of the lower semi-continuous function (see \cite{ABolte09,BPauwels16,BSTeboullle14} for more details). Let $\iota>0$ and $\Psi_{\iota}$ be the class of functions $\psi:[0,\iota)\to\mathcal{R}_+$ that satisfy the following conditions:
\begin{itemize}
	\item[(a)] $\psi(0)=0$; 
	\item[(b)] $\psi$ is positive, concave and continuous;
	\item[(c)] $\psi$ is continuously differentiable on $(0,\iota)$ with $\psi'(x)>0$ for any $x\in(0,\iota)$.
\end{itemize}

Let $g:\mathcal{R}^n\to (-\infty,\infty]$ be a given proper lower semicontinuous function. Suppose that $x\in {\rm dom}\,g:=\{x\in \Re^n\mid g(x)<\infty\}$. The Fr\'{e}chet subdifferential of $g$ at $x$ is defined as
\[
\widehat{\partial}g(x):=\left\{ h\in \mathcal{R}^n \mid \limsup_{x\neq y\to x}\frac{g(y)-g(x)-h^T(y-x)}{\|y-x\|}\ge 0 \right\}
\] 
and the limiting subdifferential, or simply the subdifferential of $g$ at $x$, is defined by 
\[
\partial g(x):=\left\{ h\in \mathcal{R}^n \mid \exists\, \{x^k\}\to x\ {\rm and}\ \{h^k\}\to h\ {\rm satisfying}\ h^k\in \widehat{\partial}g(x^k)\ \forall\, k\right\}.
\]

\begin{definition}[KL property]
	The given proper lower semicontinuous function $g:\mathcal{R}^n\to (-\infty,\infty]$ is said to have the KL property at $\bar{x}\in {\rm dom}\,g$ if there exist $\iota>0$, a neighborhood ${\cal U}$ of $\bar{x}$ and a concave function $\psi\in\Psi_{\iota}$ such that 
	\[
	\psi'(g(x)-g(\bar{x})){\rm dist}(0,\partial g(x))\ge 1\quad \forall\, x\in{\cal U}\  {\rm and}\  g(\bar{x})<g(x)<g(\bar{x})+\iota,
	\]
	where ${\rm dist}(x,Z)=\displaystyle{\min_{z\in Z}\|y-x\|}$ is the distance from a point $x$ to a nonempty closed set $Z$. The function $g$ is said to be a KL function if it has the KL property at each point of ${\rm dom}\,g$.
\end{definition}

One most  {frequently} used functions which have the KL property are the semialgebraic functions.

\begin{definition}[Semialgebraic sets and functions]
	A set in $\mathcal{R}^{n}$ is {\it semialgebraic} if it is a finite union of sets of the form
\[
\left\{ x\in\mathcal{R}^{n}\,|\, p_{i}(x)>0,\ q_{j}(x)=0,\quad i=1,\ldots,a,\ j=1,\ldots,b \right\}\,,
\]
where $p_{i}:\mathcal{R}^{n}\to\mathcal{R}$, $i=1,\ldots,a$ and $q_{j}:\mathcal{R}^{n}\to\mathcal{R}$, $j=1,\ldots,b$ are polynomials. A mapping is semialgebraic if its graph is  semialgebraic.
\end{definition}

For this class of function, we have the following useful result (cf. \cite{BDLewis07,BDLShiota07}).

\begin{proposition}\label{prop:semialgebratic-KL}
	Suppose a proper lower semicontinuous function $g:\mathcal{R}^n\to (-\infty,\infty]$ is semialgebraic, then $g$ is a KL function.
\end{proposition}

Now, we are ready to establish the global convergence of Algorithm \ref{alg:DCA} by employing a refined  global convergence result for the proximal DCA solving the DC programming with the nonsmooth DC function, which is recently developed by Liu et al. \cite{LPTakeda19}.

\begin{theorem}
Suppose that Assumption \ref{ass:bounded} holds. Let $\rho>0$ be given and $\sigma\le 1/\|\overline{C}+\rho I\|$. Suppose that $\{Y^k\}$ is the infinite sequence generated by Algorithm \ref{alg:DCA}.  Then $\{Y^k\}$ converges to a stationary point of problem \eqref{eq:QAP_SDP_p}.
\end{theorem}
\begin{proof}
	It is easy to verify that the set $\Omega\subseteq {\cal S}^q$ defined in \eqref{eq:def-Omega} is semialgebraic. Moreover, since the conjugate function  $p^*(Y):=\displaystyle\sup_{Z\in{\cal S}^q}\left\{\langle Y,Z\rangle - \|Z\|_2\right\}$ coincides with the indicator function of the unit ball of the nuclear  norm $\|\cdot\|_*$, i.e., $\{ Y\in{\cal S}^n \mid \|Y\|_*\le 1 \}$ (cf. \cite[Theorems 13.5 \& 13.2]{Rockafellar70}), we know that for the given $\sigma>0$ the corresponding  auxiliary major  function $E(Y,Z,W):=\langle \overline{C},Y\rangle +\rho\langle I, Y \rangle +\delta_{\Omega}(Y)-\langle Y,Z\rangle +p^*(Z)+\frac{1}{2\sigma}\|Y-W\|^2$, $Y,Z,W\in {\cal S}^n$ defined in \cite[(7)]{LPTakeda19} is semialgebraic. It then follows from Proposition \ref{prop:semialgebratic-KL} that $E$ is a KL function. Thus, the desired result follows from \cite[Theorem 3.1]{LPTakeda19} directly.  \qed
\end{proof}

Finally, we will show that  if the parameter $\rho>0$ is large enough, then the  sequence $\{Y^k\}$ obtained by Algorithm \ref{alg:DCA} will satisfy the the rank constraint of \eqref{eq:QAP_SDP_rank_simp} when $k$ sufficiently large.  

\begin{proposition}\label{prop:rank-feasible-DCA}
Suppose that Assumptions \ref{ass:bounded} and \ref{ass:Slater} hold. For each $k$, choose $W^k=U^k_1(U^k_1)^T\in \partial\,p(Y^{k})$, where $U^k_1\in{\cal R}^q$ is the orthonormal eigenvector with respect to the largest eigenvalue $\lambda_1(Y^k)$ of $Y^k$.  Let $\{Y^k\}$ be the sequence generated by Algorithm \ref{alg:DCA}.   
Then, there exists $\widehat{\rho}>0$ such that for any $\rho>\widehat{\rho}$ and each $k$ sufficiently large,
	\[
	{\rm rank}(Y^{k+1})\le 1,
	\]
	which implies that $Y^{k+1}$ is a feasible solution  of \eqref{eq:QAP_SDP_rank_simp}.
\end{proposition}
\begin{proof}
	For each $k$, since problem \eqref{eq:prox_mod} is convex, we know that $Y^{k+1}$ is the optimal solution of \eqref{eq:prox_mod} if and only if there exists $(y^{k+1},S^{k+1},Z^{k+1})\in \mathcal{R}^m\times {\cal S}^{q}\times {\cal S}^{q}$ such that $(Y^{k+1},y^{k+1},S^{k+1},Z^{k+1})$ satisfies the following KKT system:
	\begin{equation}\label{eq:KKT-subop}
	\left\{
	\begin{array}{l}
	\overline{C}+\rho(I-W^k)+{\cal A}^*y+S+Z+\displaystyle\frac{1}{\sigma}(Y-Y^k)=0, \\ [3pt]
	{\cal A}(Y)=b, \\ [3pt]
	{\cal S}^{q}_+\ni Y \perp S\in {\cal S}^{q}_-,\quad {\cal N}^q\ni Y \perp Z\in -{\cal N}^q.
	\end{array}
	\right.
	\end{equation}
	By the first equation of \eqref{eq:KKT-subop}, we know that for each $k$,
	\[
	S^{k+1}=-\rho(I-W^k)+M^{k+1},
	\]
	where $M^{k+1}:=-\overline{C}-{\cal A}^*y^{k+1}-Z^{k+1}-\displaystyle\frac{1}{\sigma}(Y^{k+1}-Y^k)$. By Weyl's eigenvalue inequality (see \cite{Weyl12} or \cite[Theorem 4.3.7]{HJohnson85}), we have for each $k$,
	\begin{equation}\label{eq:eig-ineqS}
	\lambda_{2}(S^{k+1})\le \lambda_2(-\rho(I-W^k))+\lambda_1(M^{k+1})=-\rho+\lambda_1(M^{k+1}),
	\end{equation}
	where the equality holds due to the fact that the eigenvalues $\lambda(-\rho(I-W^k))=(0,-\rho,\ldots,-\rho)\in {\cal R}^q$. Moreover, since for each $k$, $Y^{k}\in \Omega$ is bounded, we know that there exists a constant $\zeta>0$ such that for each $k$, $\|Y^k\|_2\le \zeta$.
	It follows from Assumption \ref{ass:Slater} that the level set of the dual problem \eqref{eq:prox_mod-dual} is a closed and bounded convex set (cf. \cite[Theorems 17 \& 18]{Rockafellar74}). Thus, we know that there exists a finite constant $\eta$ such that for $k$ sufficiently large,  $\lambda_1(-\overline{C}-{\cal A}^*y^{k+1}-Z^{k+1})\le \eta$, we have there exists a constant $\zeta>0$ such that for $k$ sufficiently large,
	\begin{eqnarray*}
		\lambda_1(M^{k+1})&\le& \lambda_1(-\overline{C}-{\cal A}^*y^{k+1}-Z^{k+1}) +\displaystyle\frac{1}{\sigma}\lambda_1\left(Y^{k+1}-Y^k\right) \\ [3pt]
		&\le& \lambda_1(-\overline{C}-{\cal A}^*y^{k+1}-Z^{k+1}) + \displaystyle\frac{1}{\sigma}\left\| Y^{k+1}-Y^k \right\|_2 \le  \eta+\frac{\zeta}{\sigma}.
	\end{eqnarray*}
	Therefore, we know from \eqref{eq:eig-ineqS} that if $\rho >\widehat{\rho}:=\max\{\eta,0\}+\frac{\zeta}{\sigma}>0$, then for $k$ sufficiently large,
	\begin{equation}\label{eq:eig-ineqS2}
	\lambda_{2}(S^{k+1})\le -\rho+\lambda_1(M^{k+1})\le -\rho+\eta+\frac{\zeta}{\sigma} <0.
	\end{equation}
	Finally, since ${\cal S}^{q}_+\ni Y^{k+1} \perp S^{k+1}\in {\cal S}^{q}_-$, by \eqref{eq:eig-ineqS2}, we obtain that for $k$ sufficiently large,
	\[
	{\rm rank}(S^{k+1})\ge q-1 \quad {\rm and} \quad {\rm rank}(Y^{k+1})+{\rm rank}(S^{k+1})\le q,
	\]
	which implies that ${\rm rank}(Y^{k+1})\le q-{\rm rank}(S^{k+1})\le 1$.  \qed
\end{proof}

\section{Numerical results}
In this section, we present numerical results for the
relaxation problem \eqref{eq:QAP_SDP_rank} solving by Algorithm \ref{alg:DCA}. All the data from QAPLIB \cite{HAQAPLIB} and `dre' instances \cite{DreznerHT05} are tested on a Window 10 workstation (6 core, Intel Xeon E5-2650 v3 @ 2.30~GH{\footnotesize{Z}}, 128 GB RAM). The size of most QAPs ranges from 12 to 60. During our experiments, SDPNAL+ version 1.0 \cite{STYZ19} is used as doubly nonnegative solver for solving the subproblems \eqref{eq:prox_mod}.  Algorithm \ref{alg:DCA} is implemented in the MATLAB 2015a platform. We measure the performance of Algorithm \ref{alg:DCA} by  
$$ \textrm{gap} := \frac{\textrm{PDCA}-\textrm{opt}}{\textrm{opt}} \times 100\% , $$
where `opt' denotes the optimal value (or best-known feasible solution) of the instance from QAPLIB, `PDCA' denotes the optimal value of the subproblem \eqref{eq:prox_mod}. 

\subsection{Penalty parameter}
\begin{figure}[htbp]
\centering
\subfigure[chr18a.]{
\includegraphics[width=5.6cm]{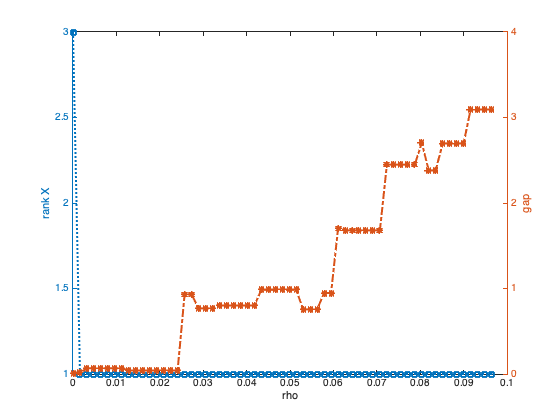}
}
\,
\subfigure[els19.]{
\includegraphics[width=5.6cm]{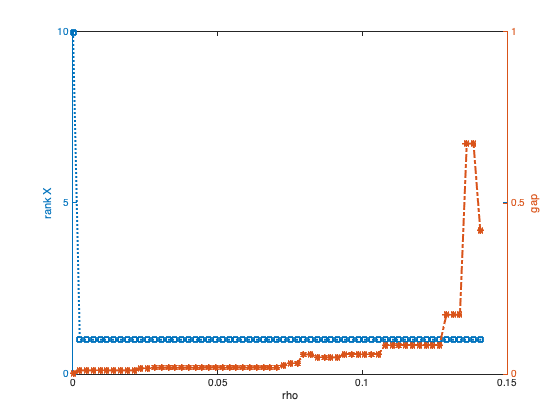}
}
\,
\subfigure[had20.]{
\includegraphics[width=5.6cm]{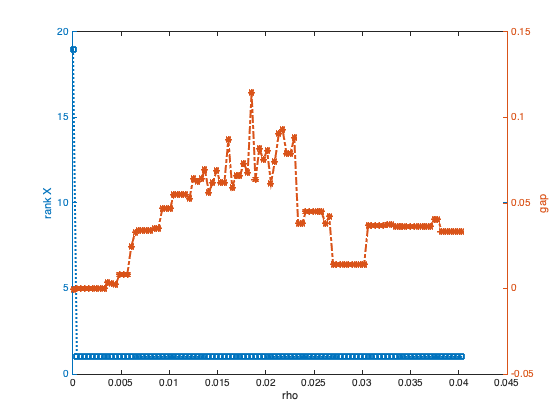}
}
\,
\subfigure[lipa30a.]{
\includegraphics[width=5.6cm]{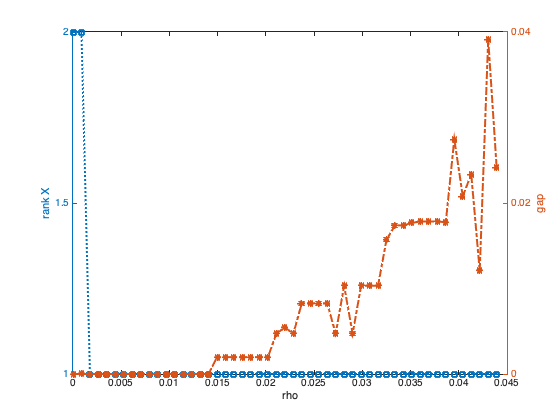}
}
\vspace{-.2 cm}
\caption{Effects of paramenters $\rho$ on gaps and ranks of solutions }
\label{fig:rho4}
\end{figure}

The penalty parameter $\rho$ is an important factor for the whole procedure of Algorithm \ref{alg:DCA}. Figure \ref{fig:rho4} shows the effect of the paramenter $\rho$ on the gaps and the ranks of the sequences generated by Algorithm \ref{alg:DCA} for chr18a, els19, had20 and lipa30a. In each subfigure, x-axis is the range of the parameter $\rho$, the left and right y-axis denote the ranks of the generated solutions and the gaps of the optimality for the different $\rho$ respectively. As shown in Fig. \ref{fig:rho4} (a) and (b), if $\rho$  increases from $0$, chr18 and els19 problems can obtain the optimal solutions of the problem \eqref{eq:QAP_SDP_rank_simp} since the gaps are zeroes.

Although larger $\rho$ can help the solutions satisfying the rank-one constraint in the problem \eqref{eq:QAP_SDP_rank_simp} (Proposition \ref{prop:rank-feasible-DCA}), the parameter $\rho$ should not be too large. In fact, as demonstrated by (c) had20 and (d) lipa30a in Fig. \ref{fig:rho4}, when $\rho$ increases larger than certain value, the gaps of these two problems oscillate up and down which imply the penalty problem \eqref{eq:QAP_SDP_p} may move away from the target problem \eqref{eq:QAP_SDP_rank_simp}. In our implementation, a bisection strategy is used for finding a suitable parameter $\rho$ for Algorithm \ref{alg:DCA}.

\subsection{Numerical performance}
Table \ref{table1} summarizes the quality of the solutions obtained by our proposed DCA approach for solving the problems from QAPLIB \cite{HAQAPLIB} and `dre' instances \cite{DreznerHT05} ($107$ instances). It can be seen from Table \ref{table1} that for $69$ instances we are able to solve the problems exactly; for $32$ instances we are able to obtain a feasible solution whose gap is less than or equal $4\%$; for $6$ instances we  obtain a feasible solution whose gap is larger than $4\%$.

\vspace{-.2 cm}
\begin{center}
\begin{table} 
	\centering   
	\caption{Summary of numerical performance of  Algorithm \ref{alg:DCA}}   
	\label{table1}  
	\begin{tabular}{|l|c|c|c|c|}  
		\hline  
		\multirow{2}{*}{Problem set (No.)} & \multicolumn{3}{c|}{gap}  &  \multirow{2}{*}{ Problem } \\ \cline{2-4}
		        & 0 & $\le 4\%$&$> 4\%$ & \\   
		\hline
		drexxx(6)& 6 & 0 & 0   &  dre15, dre18, dre21, \\
		&&&& dre24, dre30,  dre42 \\ \hline
		bur26x(8)& 0 & 8 &  0  & bur26a-h\\ \hline
		chrxxx(14)& 14& 0 &  0 &  chr12x, chr15x, chr18x,  \\
		&&&& chr20x, chr22x, chr25a\\\hline
		els19(1)& 1 & 0 &  0 & els19 \\ \hline
		escxxx(14)& 11 & 1 &  2 & esc16a-j, esc32a-g\\ \hline
		hadxx(5)& 5 & 0&  0 & had12, had14-had20 \\ \hline
		kra32x(3)& 1 & 2&  0 & kra30a-b, kra32\\ \hline
		lipaxxx(10)& 10 & 0 &  0 & lipa20x,  lipa30x, lipa40x,  \\
		&&&& lipa50x, lipa60x \\ \hline
		nugxx(13)& 8 & 5 &  0 & nug12, nug14-nug22,  \\
		&&&& nug25, nug27, nug28\\ \hline
		rouxx(3)& 3 & 0 &  0 & rou12, rou15, rou20 \\ \hline
		scrxx(3)& 3 & 0 &  0 & scr12, scr15, scr20 \\ \hline
		skoxx(5)& 0 & 2 &  3 & sko42, sko56, sko64,  \\
		&&&& sko72, sko81\\ \hline
		ste36x(3)& 0 & 3 &  0 & ste36a-c\\ \hline
		taixxx(17)& 7 & 9 &  1  & tai12x, tai15x, tai17x,  tai20x, \\
		&&&& tai25x, tai30x, tai35x,  \\
		&&&& tai40x, tai50a, tai60b\\ \hline
		thoxx(2)& 0 & 2 &  0 & tho30, tho40\\ \hline
		\hline
		Total(107)& 69 & 32 &  6 & \\ \hline
	\end{tabular}
\end{table}
\end{center}

The detail numerical results of Algorithm \ref{alg:DCA} for solving the `dre' instances from \cite{DreznerHT05} and QAPLIB \cite{HAQAPLIB} are reported in Tables \ref{table2} and \ref{table3}. In the these tables, `time' column (in hours:minutes:seconds) reports the CPU time of Algorithm \ref{alg:DCA} and `permutaion/bound' column reports the feasible solution generated by solving the relaxation problem \eqref{eq:prox_mod} of the rank-1 constrained DNN problem \eqref{eq:QAP_SDP_p}.

The `dre' problem instances \cite{DreznerHT05} are based on a rectangular grid where all nonadjacent nodes have zero weight, making the value of the objective function increase steeply with just a slight change from the optimal permutation. The `dre' instances are difficult to solve, especially for many metaheuristic-based methods, since they are ill-conditioned and hard to break out the `basin' of the local minimal.  The best known solutions for the `dre' problems have been found by branch and bound in \cite{DreznerHT05}. Notably, by employing our proposed DCA based approach  Algorithm \ref{alg:DCA}, we are able to obtain the global optimal solutions of the `dre' problems quite efficiently. For instance, we are able to solve the instance `dre42' by  Algorithm \ref{alg:DCA} exactly in $13$ minutes. 

\vspace{-.4cm}
\begin{center}
\footnotesize
\begin{longtable}{|l|rrr|r|l|}
\caption{Numerical performance of the `dre' problem instances \cite{DreznerHT05}}   
	\label{table2}\\
\hline
 & & & & & \\[-4pt]
\multicolumn{1}{|c|}{Problem} & \multicolumn{1}{c}{opt}&\multicolumn{1}{r}{PDCA} &  \multicolumn{1}{c}{gap ($\%$)} & \multicolumn{1}{|c}{time} & \multicolumn{1}{|c|}{permutation$\slash$bound}\\
 & & & & & \\[-4pt]
\hline
\endfirsthead
\multicolumn{6}{c}{\tablename\ \thetable\ -- continued from previous page} \\
\hline

 & & & & & \\[-4pt]
\multicolumn{1}{|c|}{Problem} & \multicolumn{1}{r}{opt}&\multicolumn{1}{r}{PDCA} &  \multicolumn{1}{c}{gap($\%$)} & \multicolumn{1}{|c}{time} & \multicolumn{1}{|c|}{permutation$\slash$bound}\\
 & & & & & \\[-4pt]
\hline
\endhead

\hline  \multicolumn{6}{|c|}{{Continued on next page}} \\ \hline
\endfoot

\hline \hline
\endlastfoot

dre15   &  306  &  306  &  0 &11 &  1 13     4     6     7     9    11     5    12    14     1    15    10 2     3     8  \\ \hline 
 
dre18   &  332  &  332  &  0 &16 & 4    14    18     9    10    12     2    15     7     3     5     8     6    11    13    17      \\[1pt]  
 & & & &   &    1    16  \\ \hline 
 
 dre21  &  356  &  356  &  0 & 33 &  5     8    17    18    12    13     1    11     3     9    16     4     6    20     7    19        \\[1pt]  
 & & & &   &    14   10    15     2    21  \\ \hline 

 dre24  &  396  &  396  &  0 & 15&  3    23    14    21    22    10    16     9     7     5     8    18    13     4     2    17        \\[1pt]  
 & & & &   &     1   19    12    11    15    24     6    20 \\ \hline 

 dre30  &  508  &  508  &  0 & 1:22 &  28     2     1    17     6     3    11    21    19    22    24     8    26    20    23    13    \\[1pt]  
 & & & &   &    4    29    18    25    10    30    16    15    14    12     7     5    27     9 \\ \hline 

   dre42 &  764  &  764  &  0 & 13:00 &  3    36    41    28    30    14    34    32    42    37    33    10    27    12    35     9    \\[1pt]  
 & & & &   &    7   21     5    29    18    11     8    38    24     2    15    22     6     1    13    19  \\[1pt]  
 & & & &   &   40    23   25    39    31    16    17    26     4    20
  \\ \hline

\end{longtable}

\end{center}
\vspace{-.4cm}
In Table \ref{table3},  the upper bounds generated by Algorithm \ref{alg:DCA} are compared with the state of the art optimal values (or the best known upper bounds) in QAPLIB. Except bur$xxx$ and sko$xx$ cases, we find that most instances can either be solved exactly or achieve an upper bound which is accurate up to a relative error of $5\%$ through the penalized DC relaxation. Because the subproblems of the corresponding penalized DC problems are failed to achieve the stopping criteria $10^{-6}$ of SDPNAL+, Algorithm \ref{alg:DCA} only provides the feasible solutions for bur$xxx$ cases.  We note that the QAPLIB bounds were typically achieved using a rather large collection of different algorithms, which generally involve a branch and bound procedure requiring multiple convex relaxations, while our results are achieved by using a single relaxation. 
 
\begin{center}
\footnotesize
\begin{longtable}{|l|rrr|r|l|}
\caption{Numerical performance of the QAPLIB instances}   
	\label{table3}\\
\hline
 & & & & & \\[-4pt]
\multicolumn{1}{|c|}{Problem} & \multicolumn{1}{c}{opt}&\multicolumn{1}{r}{PDCA} &  \multicolumn{1}{c}{gap ($\%$)} & \multicolumn{1}{|c}{time} & \multicolumn{1}{|c|}{permutation$\slash$bound}\\
 & & & & & \\[-4pt]
\hline
\endfirsthead
\multicolumn{6}{c}{\tablename\ \thetable\ -- continued from previous page} \\
\hline

 & & & & & \\[-4pt]
\multicolumn{1}{|c|}{Problem} & \multicolumn{1}{r}{opt}&\multicolumn{1}{r}{PDCA} &  \multicolumn{1}{c}{gap($\%$)} & \multicolumn{1}{|c}{time} & \multicolumn{1}{|c|}{permutation$\slash$bound}\\
 & & & & & \\[-4pt]
\hline
\endhead

\hline  \multicolumn{6}{|c|}{{Continued on next page}} \\ \hline
\endfoot

\hline \hline
\endlastfoot

  bur26a   &  5426670  &  5566175  &  2.57 & 2:33 & 11  26  15  7  4  13  12  6  2  18  1  9  5  21  8  14 \\[1pt] 
 & & & &   &  3   19   20   17   10   25   24   16   23   22 \\ \hline 
 
  bur26b   &  3817852  &  3956961  &  3.64 & 1:55 & 15  16  10  7  4  12  2  23  22  18  5  9  1  21  8  14 \\[1pt] 
 & & & &   &  3   20   19   25   11   26   24   17   6   13 \\ \hline 
 
  bur26c   &  5426795  &  5523812  &  1.79 & 7:03 & 13  3  12  7  16  11  25  10  15  9  8  19  18  20  4  21 \\[1pt] 
 & & & &   &  1   14   5   6   22   24   2   23   26   17 \\ \hline 
 
  bur26d   &  3821225  &  3902248  &  2.12 & 3:59 & 22  23  3  2  16  11  17  21  15  9  8  18  19  20  12  25 \\[1pt] 
 & & & &   &  14   1   5   13   24   6   4   7   26   10 \\ \hline 
 
  bur26e   &  5386879  &  5470899  &  1.56 & 6:26 & 22  3  6  7  12  26  1  16  11  9  18  19  20  14  13  8 \\[1pt] 
 & & & &   &  5   15   21   2   17   24   4   10   25   23 \\ \hline 
 
  bur26f   &  3782044  &  3847551  &  1.73 & 6:16 & 6  22  4  3  12  25  7  1  23  15  20  18  19  14  16  10 \\[1pt] 
 & & & &   &  5   21   9   24   2   17   26   13   11   8 \\ \hline 
 
  bur26g   &  10117172  &  10332466  &  2.13 & 7:32 & 2  11  22  23  13  10  25  8  1  21  20  4  7  18  12  15 \\[1pt] 
 & & & &   &  9   19   5   26   16   6   14   3   24   17 \\ \hline 
 
  bur26h   &  7098658  &  7257159  &  2.23 & 4:06 & 22  16  13  26  14  10  21  1  8  15  4  20  18  7  12  17 \\[1pt] 
 & & & &   &  19   5   9   2   11   3   23   6   24   25 \\ \hline 
 
  chr12a   &  9552  &  9552  &  0 & 04 & 7  5  12  2  1  3  9  11  10  6  8  4 \\ \hline 
 
  chr12b   &  9742  &  9742  &  0 & 04 & 5  7  1  10  11  3  4  2  9  6  12  8 \\ \hline 
 
  chr12c   &  11156  &  11156  &  0 & 04 & 7  5  1  3  10  4  8  6  9  11  2  12 \\ \hline 
 
  chr15a   &  9896  &  9896  &  0 & 07 & 5  10  8  13  12  11  14  2  4  6  7  15  3  1  9 \\ \hline 
 
  chr15b   &  7990  &  7990  &  0 & 07 & 4  13  15  1  9  2  5  12  6  14  7  3  10  11  8 \\ \hline 
 
  chr15c   &  9504  &  9504  &  0 & 07 & 13  2  5  7  8  1  14  6  4  3  15  9  12  11  10 \\ \hline 
 
  chr18a   &  11098  &  11098  &  0 & 13 & 3  13  6  4  18  12  10  5  1  11  8  7  17  14  9  16 \\[1pt] 
 & & & &   &  15   2 \\ \hline 
 
  chr18b   &  1534  &  1534  &  0 & 1:52 & 10  7  11  8  12  9  15  6  14  5  13  4  16  1  17  2 \\[1pt] 
 & & & &   &  18   3 \\ \hline 
 
  chr20a   &  2192  &  2192  &  0 & 27 & 3  20  7  18  9  12  19  4  10  11  1  6  15  8  2  5 \\[1pt] 
 & & & &   &  14   16   13   17 \\ \hline 
 
  chr20b   &  2298  &  2298  &  0 & 28 & 20  3  9  7  1  12  16  6  8  14  10  4  5  13  17  2 \\[1pt] 
 & & & &   &  18   11   19   15 \\ \hline 
 
  chr20c   &  14142  &  14142  &  0 & 28 & 12  6  9  2  10  11  3  4  15  18  7  13  16  5  14  17 \\[1pt] 
 & & & &   &  19   1   8   20 \\ \hline 
 
  chr22a   &  6156  &  6156  &  0 & 39 & 15  2  21  8  16  1  7  18  14  13  5  17  6  11  3  4 \\[1pt] 
 & & & &   &  20   19   9   22   10   12 \\ \hline 
 
  chr22b   &  6194  &  6194  &  0 & 28 & 10  19  3  1  20  2  6  4  7  8  17  13  11  15  21  12 \\[1pt] 
 & & & &   &  9   5   22   14   18   16 \\ \hline 
 
  chr25a   &  3796  &  3796  &  0 & 2:00 & 25  12  5  3  18  4  16  8  20  10  14  6  23  15  24  19 \\[1pt] 
 & & & &   &  13   1   21   11   17   2   22   7   9 \\ \hline 
 
  els19   &  17212548  &  17212548  &  0 & 31 & 9  10  7  19  14  18  13  17  6  11  4  5  12  8  16  15 \\[1pt] 
 & & & &   &  1   2   3 \\ \hline 
 
  esc16a   &  68  &  68  &  0 &  18  & 12  15  7  11  14  6  10  8  4  16  13  3  5  9  2  1 \\ \hline 
 
  esc16b   &  292  &  292  &  0 & 40 & 7  6  8  14  16  10  13  2  15  4  12  5  1  11  9  3 \\ \hline 
 
  esc16c   &  160  &  160  &  0 & 1:00 & 15  10  9  2  1  14  11  3  4  13  8  16  12  6  7  5 \\ \hline 
 
  esc16d   &  16  &  16  &  0 & 8 & 14  6  7  10  13  5  16  2  4  1  3  12  15  11  8  9 \\ \hline 
 
  esc16e   &  28  &  30  &  7.14 &  15 & 15  10  7  14  4  8  6  16  12  2  1  5  9  13  11  3 \\ \hline 
 
  esc16g   &  26  &  26  &  0 & 13 & 7  12  10  16  4  8  6  14  11  3  2  1  15  13  5  9 \\ \hline 
 
  esc16h   &  996  &  996  &  0 &  18  & 6  5  13  7  12  11  15  4  8  3  16  9  2  1  10  14 \\ \hline 
 
  esc16i   &  14  &  14  &  0 &  19 & 7  5  3  1  9  10  12  2  4  6  11  13  15  8  14  16 \\ \hline 
 
  esc16j   &  8  &  8  &  0 & 11 & 11  4  5  8  14  16  13  9  7  1  10  12  15  3  6  2 \\ \hline 
 
  esc32a   &  130  &  150  &  15.38 & 13:11 & 28  12  27  19  4  18  16  21  11  32  14  8  26  10  25  23 \\[1pt] 
 & & & &   &  6   9   13   17   15   22   7   31   30   24   29   5   3   1   20   2 \\ \hline 
 
  esc32b   &  168  &  168  &  0 & 12:18 & 1  3  9  25  11  27  5  7  6  2  8  4  13  29  14  10 \\[1pt] 
 & & & &   &  30   26   15   31   16   12   32   28   21   19   23   20   18   17   22   24 \\ \hline 
 
  esc32c   &  642  &  646  &  0.62 &    2:26:27  &  11    30    14    15    17    24    32     1    23     4     2    18    22     6    21     5   \\[1pt] 
 & & & &   &   8     7    20    19    29     9     3    13    26    12    16    31    25    27    10   28\\ \hline 
 
%
  esc32e   &  2  &  2  &  0 &  12:36 & 19  2  17  9  3  7  12  30  16  10  14  13  31  15  21  18 \\[1pt] 
 & & & &   &  22   28   4   26   6   25   24   11   20   8   27   5   1   32   23   29 \\ \hline 
 
  esc32g   &  6  &  6  &  0 & 12:23 & 2  13  5  6  14  22  20  15  16  26  28  8  32  23  31  1 \\[1pt] 
 & & & &   &  11   17   12   24   3   7   19   4   30   25   21   10   27   29   18   9 \\ \hline 
 
  had12   &  1652  &  1652  &  0 & 4 & 3  10  11  2  12  5  7  6  8  1  4  9 \\ \hline 
 
  had14   &  2724  &  2724  &  0 & 5 & 8  13  10  5  12  11  2  14  3  6  7  1  9  4 \\ \hline 
 
  had16   &  3720  &  3720  &  0 &  8 & 9  4  16  1  7  8  6  14  15  11  12  10  5  3  2  13 \\ \hline 
 
  had18   &  5358  &  5358  &  0 & 11 & 8  15  16  14  7  18  6  11  1  10  12  5  3  13  2  17 \\[1pt] 
 & & & &   &  9   4 \\ \hline 
 
  had20   &  6922  &  6922  &  0 &  42 & 8  15  1  14  19  6  7  17  16  12  10  11  5  20  2  3 \\[1pt] 
 & & & &   &  4   9   18   13 \\ \hline 
 
  kra30a   &  88900  &  88900  &  0 &  23:16 & 9  13  28  27  8  7  10  30  20  21  23  19  24  29  14  1 \\[1pt] 
 & & & &   &  11   12   18   16   17   22   26   2   5   3   25   6   4   15 \\ \hline 
 
  kra30b   &  91420  &  92010  &  0.65 & 14:45 & 24  15  22  6  4  3  5  2  17  1  16  11  25  26  30  19 \\[1pt] 
 & & & &   &  10   28   29   21   23   20   8   12   9   7   13   27   18   14 \\ \hline 
 
  kra32   &  88900  &  89100  &  0.22 &  2:11:24&  8    10     5     9     6    14    27    23    24    31    12    16    28     4    13      7  \\[1pt] 
 & & & &   &  11     3    29    22    26    15    30    21     1     2    20    32    25    18    19   17 \\ \hline 
 
  lipa20a   &  3683  &  3683  &  0 & 18  & 19  17  7  1  5  9  10  12  4  16  20  6  3  14  11  15 \\[1pt] 
 & & & &   &  13   8   2   18 \\ \hline 
 
  lipa20b   &  27076  &  27076  &  0 & 19 & 1  2  3  4  5  6  7  8  9  10  11  12  13  14  15  16 \\[1pt] 
 & & & &   &  17   18   19   20 \\ \hline 
 
  lipa30a   &  13178  &  13178  &  0 & 4:43 & 9  13  22  17  25  23  29  12  11  6  5  28  20  27  14  4 \\[1pt] 
 & & & &   &  18   8   19   30   21   7   15   24   26   3   16   1   10   2 \\ \hline 
 
  lipa30b   &  151426  &  151426  &  0 &  4:51 & 1  2  3  4  5  6  7  8  9  10  11  12  13  14  15  16 \\[1pt] 
 & & & &   &  17   18   19   20   21   22   23   24   25   26   27   28   29   30 \\ \hline 
 
  lipa40a   &  31538  &  31538  &  0 & 34:01  & 7  6  14  27  19  37  21  2  16  1  40  24  30  23  5  28 \\[1pt] 
 & & & &   & 22  8  20  35  32  26  29  3  4  11  36  10  13  38  9  17 \\[1pt] 
 & & & &   &  31   18   33   15   25   34   39   12 \\ \hline 
 
  lipa40b   &  476581  &  476581  &  0 &  40:58 & 1  2  3  4  5  6  7  8  9  10  11  12  13  14  15  16 \\[1pt] 
 & & & &   & 17  18  19  20  21  22  23  24  25  26  27  28  29  30  31  32 \\[1pt] 
 & & & &   &  33   34   35   36   37   38   39   40 \\ \hline 
 
  lipa50a   &  62093  &  62093  &  0 & 20:30 & 28  32  37  39  49  23  19  44  33  7  14  30  15  5  36  6 \\[1pt] 
 & & & &   & 17  26  48  25  40  3  45  27  18  31  29  16  9  12  1  8 \\[1pt] 
 & & & &   & 4  2  50  21  43  35  24  38  34  46  42  13  20  22  41  47 \\[1pt] 
 & & & &   &  10   11 \\ \hline 
 
  lipa50b   &  1210244  &  1210244  &  0 & 23:04 & 1  2  3  4  5  6  7  8  9  10  11  12  13  14  15  16 \\[1pt] 
 & & & &   & 17  18  19  20  21  22  23  24  25  26  27  28  29  30  31  32 \\[1pt] 
 & & & &   & 33  34  35  36  37  38  39  40  41  42  43  44  45  46  47  48 \\[1pt] 
 & & & &   &  49   50 \\ \hline 
 
  lipa60a   &  107218  &  107218  &  0 & 2:44:09 & 25  48  17  4  50  13  16  41  1  37  22  27  46  34  38  12 \\[1pt] 
 & & & &   & 9  3  8  33  47  54  31  43  40  10  35  23  29  57  2  6 \\[1pt] 
 & & & &   & 51  56  49  21  30  36  15  39  59  18  52  28  26  44  14  60 \\[1pt] 
 & & & &   &  32   5   20   11   55   45   53   24   42   7   58   19 \\ \hline 
 
  lipa60b   &  2520135  &  2520135  &  0 & 36:07 & 1  2  3  4  5  6  7  8  9  10  11  12  13  14  15  16 \\[1pt] 
 & & & &   & 17  18  19  20  21  22  23  24  25  26  27  28  29  30  31  32 \\[1pt] 
 & & & &   & 33  34  35  36  37  38  39  40  41  42  43  44  45  46  47  48 \\[1pt] 
 & & & &   &  49   50   51   52   53   54   55   56   57   58   59   60 \\ \hline 
 
  nug12   &  578  &  578  &  0 &  12 & 3  9  7  12  1  11  8  4  2  10  6  5 \\ \hline 
 
  nug14   &  1014  &  1014  &  0  & 4:40 & 9     8    13     2     1    11     7    14     3     4    12     5     6    10 \\ \hline 
 
  nug15   &  1150  &  1150  &  0 & 5:19 & 12  5  6  15  10  11  7  14  3  4  9  8  13  2  1 \\ \hline 
 
  nug16a   &  1610  &  1610  &  0 &  32 & 9  14  2  15  16  3  10  12  8  11  6  5  7  1  4  13 \\ \hline 
 
  nug16b   &  1240  &  1240  &  0 & 9:27 & 8  11  3  5  13  9  7  1  12  2  10  6  16  4  15  14 \\ \hline 
 
  nug17   &  1732  &  1732  &  0 &  10:17 &16    15     2    14     9    11     8    12    10     3     4     1     7     6    13    17 \\[1pt] 
 & & & &   &  5 \\ \hline 
 
  nug18   &  1930  &  1936  &  0.31 &  29 & 17  1  7  5  6  11  8  12  10  13  4  15  2  3  9  16 \\[1pt] 
 & & & &   &  18   14 \\ \hline 
 
  nug20   &  2570  &  2570  &  0 & 2:40 & 18  14  10  3  9  4  2  12  11  16  19  15  20  8  13  17 \\[1pt] 
 & & & &   &  5   7   1   6 \\ \hline 
 
  nug21   &  2438  &  2442  &  0.16 &  17:25 & 4  18  11  16  13  6  15  14  19  8  7  1  12  17  20  21 \\[1pt] 
 & & & &   &  3   9   10   2   5 \\ \hline 
 
  nug22   &  3596  &  3642  &  1.28 & 4:02 & 20  4  16  8  19  7  10  9  13  21  17  14  15  11  3  18 \\[1pt] 
 & & & &   &  22   1   12   6   2   5 \\ \hline 
 
  nug25   &  3744  &  3750  &  0.16 & 35:06 & 12  4  14  23  13  24  21  7  10  20  17  16  6  19  1  11 \\[1pt] 
 & & & &   &  25   8   9   15   18   2   3   22   5 \\ \hline 
 
  nug27   &  5234  &  5236  &  0.04 & 30:21 & 23  18  3  1  27  17  5  13  7  15  12  26  8  19  20  2 \\[1pt] 
 & & & &   &  24   21   4   10   9   14   22   25   6   16   11 \\ \hline 
 
  nug28   &  5166  &  5166  &  0 & 35:33 & 11     8    20    28     1     9    18    26    16    17    19    10    15     7    14    27 \\[1pt] 
 & & & &   &   4    13    25     6    22    12     5     3    24     2    21    23 \\ \hline  

  rou12   &  235528  &  235528  &  0 & 04 & 6  5  11  9  2  8  3  1  12  7  4  10 \\ \hline 
 
  rou15   &  354210  &  354210  &  0 & 5:46 & 12  6  8  13  5  3  15  2  7  1  9  10  4  14  11 \\ \hline 
 
  rou20   &  725522  &  725582  &  0 & 4:54 & 12  11  9  17  7  14  8  13  6  10  19  18  5  2  16  3 \\[1pt] 
 & & & &   &  4   1   20   15 \\ \hline 
 
  scr12   &  31410  &  31410  &  0 & 04 & 8  6  11  10  2  9  5  1  12  7  4  3 \\ \hline 
 
  scr15   &  51140  &  51140  &  0 & 06 & 12  10  11  14  1  13  9  5  15  6  4  2  8  7  3 \\ \hline 
 
  scr20   &  110030  &  110030  &  0 & 13:42 & 20  7  12  6  4  8  3  2  14  11  18  9  19  15  16  17 \\[1pt] 
 & & & &   &  13   5   10   1 \\ \hline 
 
  sko42   &  15812  &  15952  &  0.89 &   3:58:53 & 23  16  30  36  3  21  8  22  39  37  24  1  38  6  28  2 \\[1pt] 
 & & & &   & 14  18  32  40  25  33  4  35  26  9  15  29  27  10  19  20 \\[1pt] 
 & & & &   &  34   42   11   13   31   5   12   7   17   41 \\ \hline

  sko56   &  34458  &  36490  & 5.90 &   3:18:03 &  21    42     2    11    30    51    49     3    25     6    56    33     9    23    13    54  \\[1pt] 
 & & & &   &  38    24     8    50    32    36    44    10    31    45    53    37    27    12    14    41  \\[1pt] 
 & & & &   &   43    55    22     5    29    20    52    35    26    46     4    28    34    48    18    40    \\[1pt] 
 & & & &   &       19    15    16     1    47    17     7    39 \\ \hline

 sko64   &  48498  &  50948  & 5.05 &   3:48:56 &  31    57    56    53    51     3    43    15     2    54    17    52    23    32    11    27  \\[1pt] 
 & & & &   &     63    26    59    49    13     6    40    46    33     7    44    16    64    39    38    61      \\[1pt] 
 & & & &   &      10    42    45    24     4    62    30    41    35    37     1    21    55     5    28    34     \\[1pt] 
 & & & &   &       9    25    58    22    19     8    48    12    18    50    14    20    29    47    60    36
 \\ \hline 

 sko72   &  66256  &  70318  & 6.13 &   14:42:49 &  45    11    51    71    29    14    60    31    12    53    41    50    57    44    28    21    \\[1pt] 
 & & & &   &    10    64    25    52    68     1     6    13    47    67    35    43    20    16    59    49      \\[1pt] 
 & & & &   &      7    22    23    27    32     4    63    55    33    46    65    17    30    24    69    66     \\[1pt] 
 & & & &   &         36     9    38    37     2    34     3    72    42    19    39    61    26    62    56    54  \\[1pt] 
 & & & &   &      5    58    15    40    48    70     8    18  \\ \hline 
      
 sko81   &  90998  &  93356  & 2.59 &   25:47:17 &    47    67    80    65    34    39    69    74    30    40    23    63    38    20    33    26   \\[1pt] 
 & & & &   &   81    25    54     4    79  35    51     2    12    14     3    77    24    73    58    32       \\[1pt] 
 & & & &   &    70    21    78    61    10    19    75    37    18    15     49    42    72    22    43    59    \\[1pt] 
 & & & &   &      13    44     7    41     6    28    71    55    62     1    11    31     5    36     9   53 \\[1pt] 
 & & & &   &      29    46    52    27    45    56    17    50     8    60    64    76    16    68    57    48 \\[1pt] 
 & & & &   &     66  \\ \hline          
 
  ste36a   &  9526  &  9640  &  1.20 &   1:12:45 & 35  16  1  15  14  28  29  30  31  17  18  10  7  11  20  19 \\[1pt] 
 & & & &   & 32  34  2  8  4  13  12  23  22  21  33  36  3  9  5  6 \\[1pt] 
 & & & &   &  27   26   25   24 \\ \hline 
 
  ste36b   &  15852  &  15932  &  0.50 & 25:14 & 35  31  30  29  28  1  15  9  16  33  34  32  19  20  7  10 \\[1pt] 
 & & & &   & 18  17  26  25  23  14  11  13  4  8  2  24  22  21  27  12 \\[1pt] 
 & & & &   &  6   5   3   36 \\ \hline 
 
  ste36c   &  8239110  &  8254628  &  0.19 &   2:19:21 & 24  25  26  27  11  6  5  3  35  22  21  23  14  12  13  4 \\[1pt] 
 & & & &   & 8  2  33  32  19  28  20  7  10  18  17  34  31  30  29  15 \\[1pt] 
 & & & &   &  1   9   16   36 \\ \hline 
 
  tai12a   &  224416  &  224416  &  0 & 05 & 8  1  6  2  11  10  3  5  9  7  12  4 \\ \hline 
 
  tai12b   &  39464925  &  39464925  &  0 & 4 & 9  4  6  3  11  7  12  2  8  10  1  5 \\ \hline 
 
  tai15a   &  388214  &  388870  &  0.17 & 31 & 6  10  4  7  2  9  1  11  3  14  13  15  5  12  8 \\ \hline 
 
  tai15b   &  51765268  &  51765268  &  0 & 15 & 1  9  4  6  8  15  7  11  3  5  2  14  13  12  10 \\ \hline 
 
  tai17a   &  491812  &  491812  &  0 & 6:18 & 12  2  6  7  4  8  14  5  11  3  16  13  17  9  1  10 \\[1pt] 
 & & & &   &  15 \\ \hline 
 
  tai20a   &  703482  &  703482  &  0 & 14:48 & 10  9  12  20  19  3  14  6  17  11  5  7  15  16  18  2 \\[1pt] 
 & & & &   &  4   8   13   1 \\ \hline 
 
  tai20b   &  122455319  &  122455319  &  0 & 48 & 8  16  14  17  4  11  3  19  7  9  1  15  6  13  10  2 \\[1pt] 
 & & & &   &  5   20   18   12 \\ \hline 
 
  tai25a   &  1167256  &  1217842  &  4.33 & 31:47 &  20     1    10     7     9    13     4    19     3     2    11    15     8    21    12    14    \\[1pt] 
 & & & &   &  18    25    23    17     5    22    24     6    16 \\ \hline 
 
  tai25b   &  344355646  &  344855160  &  0.15 & 16:03 & 4  25  16  9  13  5  6  19  7  17  10  3  15  20  18  2 \\[1pt] 
 & & & &   &  22   23   8   11   21   24   14   12   1 \\ \hline 
 
  tai30a   &  1818146  &  1818146  &  0 & 33:33 & 19  18  4  24  30  25  5  7  1  22  28  20  11  13  9  16 \\[1pt] 
 & & & &   &  8   10   17   21   12   29   2   15   3   14   26   27   23   6 \\ \hline 
 
  tai30b   &  637117113  &  644555585  &  1.17 & 4:52 & 4  15  5  8  21  11  30  14  17  20  6  13  18  7  23  10 \\[1pt] 
 & & & &   &  24   27   29   9   19   28   2   26   12   22   25   16   1   3 \\ \hline 
 
  tai35a   &  2422002  &  2431214  &  0.38 & 37:01 & 35  9  21  7  23  2  12  8  20  26  33  16  18  24  22  15 \\[1pt] 
 & & & &   & 30  3  14  11  1  5  17  4  10  13  27  29  34  32  25  31 \\[1pt] 
 & & & &   &  28   6   19 \\ \hline 
 
  tai35b   &  283315445  &  284614706  &  0.46 &   1:14:45 & 1    17     5    30    11     2    22     3    14    10     9    20    18    12    19    33  \\[1pt] 
 & & & &   &   8    32    13    27    15    21    34    35     7    24     6    23    31    28    16     4  \\[1pt] 
 & & & &   &  29   26   25 \\ \hline 
  
  tai40a   &  3139370  &  3154106  &  0.47 &   1:54:41 & 6  14  15  12  9  3  19  20  28  27  5  2  7  31  36  4 \\[1pt] 
 & & & &   & 37  13  29  35  38  32  11  1  22  39  18  30  23  33  25  24 \\[1pt] 
 & & & &   &  21   8   10   17   34   16   40   26 \\ \hline

  tai40b   &  637250948  &  640933239  &  0.46 &   2:01:44 & 19     1    25    11    37    31    36    15    39    13    22     7    40    27     4    33 \\[1pt] 
 & & & &   & 14    34     9    10    38    23     5    32    35    12    16     3    24     2    21    28  \\[1pt] 
 & & & &   &   20    17    26     6    30    29    18     8 \\ \hline

  tai50a   &  4938796  &  5086610  &  2.99 &    11:54:27 & 5    42    35    34    36    37     4     7    21    10     8    25    48    27    41    11  \\[1pt] 
 & & & &   &  18    32     6    44    12    33    15     2    29    50    20     3    43    14    46    38   \\[1pt] 
 & & & &   &   26    40    47    39    19    28    45    49    13    30    23     1    16    22    17    31  \\[1pt] 
 & & & &   &    9    24 \\ \hline 
     
tai60b   &  7205962  &  7417848  & 2.94 &       5:41:01   &  18    41    38     7    20    26    30    36    16    31    59    57    21    54    19    48    \\[1pt] 
 & & & &   &   23    47    27     4    10    14     5    43    58     3     6    50    12    49     8    44     \\[1pt] 
 & & & &   &   17    29    15    22    33    46    35    32    25    52    34    39    45    55    60    11   \\[1pt]  
 & & & &   &    1    13    53    37    51     9    40    28    24     2    42    56 \\ \hline

tho30   &  149936  &  151156  &  0.81 & 26:26 & 29  14  1  2  8  28  22  25  12  20  19  11  24  27  17  26 \\[1pt] 
 & & & &   &  30   10   6   15   3   7   5   4   16   23   21   13   9   18 \\ \hline

tho40   &  240516  &  242282  &  0.73 &   5:38:41 & 38  37  5  13  26  27  35  31  4  28  9  32  21  8  29  25 \\[1pt] 
 & & & &   & 18  33  22  16  30  6  12  34  39  14  20  15  1  10  11  17 \\[1pt] 
 & & & &   &  19   2   40   23   24   7   36   3 \\ \hline

\end{longtable}
\end{center}

\section{Conclusion}\label{sec:conclusion}
This paper established an exact rank constrained DNN formulation of QAP. Under the framework of DC programming, we are able to solve the penalized DC problem efficiently by the semi-proximal augmented Lagrangian method. If the subproblems can be solved successfully, our algorithm usually reaches the optimal solutions of QAP exactly. Even if the subproblem is difficult to solve, our proposed algorithm still can provide a good feasible solution close to the optimal upper bound in QAPLIB. As a future work, we will investigate the structure of the constraints of the penalized DC problem and try to reduce the number of constraints for solving the rank constrained DNN formulation of QAPs more efficiently.

\begin{acknowledgements}
We would like to thank Dr. Xudong Li and Dr. Ying Cui for  many helpful discussions on this work.
\end{acknowledgements}

%
%



\end{document}